\documentclass{article}
\usepackage[utf8]{inputenc}
\usepackage{amsmath}
\usepackage{amssymb}
\usepackage{graphicx}
\usepackage{epstopdf}
\usepackage{amsthm}
\usepackage{hyperref}
\usepackage{amsfonts}
\usepackage{array}
\usepackage{amsthm}
\usepackage{tabularx}
\usepackage{cite}
\usepackage{url}
\usepackage{setspace}
\usepackage{geometry}
\usepackage{subfigure}
\usepackage{authblk}
\usepackage{hyperref}
\usepackage{multirow}
\usepackage{booktabs}
\usepackage[ruled,vlined]{algorithm2e}
\newtheorem{theorem}{Theorem}
\newtheorem{remark}{Remark}
\newtheorem{definition}{Definition}

\newtheorem{lemma}{Lemma}
\newtheorem{example}{Example}
\newtheorem{corollary}{Corollary}

\title{Semi-parametric estimation and prediction intervals in state space models}
\begin{document}
\author{Yunyi Zhang, Tingting Wang, Dimitris N. Politis}

\maketitle
\abstract{Literatures in state space models focus on parametric inference and prediction, which fail if the state space model is not fully specified and the maximum likelihood estimation does not work. In this paper, we assume the state transition matrix and the distribution of state noises are unknown. Under this assumption, we provide methods to consistently estimate these terms. In addition, we introduce an algorithm to construct consistent prediction intervals for state vectors and future observations. We complement the asymptotic results with several numerical experiments.}
\maketitle
\section{Introduction}
State space models (along with the Kalman filter \cite{KalmanFilter}) were originally developed to control linear systems (Davis and Vinter \cite{LinearSystem}) but have been widely used in many disciplines, including econometrics (Stock and Watson \cite{STOCK2016415} and Hamilton \cite{HAMILTON19943039}), computer science (Rangapuram et.al. \cite{NEURIPS2018_5cf68969}), engineering \cite{KRISTIANSEN20052195} and etc. For some practical applications, Scott et.al. \cite{5466306} applied Kalman filter to obtain information from video scenes, Suk.et.al. \cite{SUK2016292} applied state space model to analyze fMRI and Glickman and Stern \cite{doi:10.1080/01621459.1998.10474084} developed a model to predict NFL scores.

State space models are powerful tools for statisticians as well. For instance, in time series literature, linear state space model can be used to find the best linear mean-square predictors for the future observations (e.g. chapter 12.2, Brockwell and Davis \cite{time_series_analysis}). There are lots of researches in state space models, including maximum likelihood estimation(chapter 7 in Durbin and Koopman \cite{Likelihood} and Jong \cite{10.1093/biomet/75.1.165}) and Bayesian estimation (chapter 13 in \cite{HMM}) of parameters. Classical state space models are linear and have normal noises, but Carter and Kohn \cite{10.2307/2337125} introduced how to perform inference on non-Gaussian state space models through Gibbs sampling. Briers et.al. \cite{Filter} introduced how to smooth a non-Gaussian and non-linear state space model. Advancements in sampling techniques like Pompe et.al. \cite{pompe2020} also bring new insights in this field.

Despite great achievements in parametric and Bayesian settings, to the best of our knowledge, few researches on state space models have been conducted under non-parametric settings. However, there are interesting problems in this field. For example, by adopting notations in chapter 1.3.3, \cite{HMM}, if we do not know one or many of the following terms: state transition matrix, measurement transition matrix, distribution of state noise and distribution of measurement noise, we hope to find consistent estimators of them. Moreover, major usages of state space models include estimating state vectors based on observations, so we hope to know prediction intervals for state vectors or future observations apart from predictors in the nonparametric prediction setting introduced in Stein \cite{doi:10.1080/01621459.1985.10478220}, Thombs and Schucany \cite{10.2307/2289788}, Pan and Politis \cite{PAN20161} and Politis \cite{model_free}. Anderson et.al. \cite{10.2307/2239521} provided a consistent estimator for state transition matrix of a special state space model, Costa and Alpuim \cite{COSTA20101889} tried to estimate some parameters in the state space model non-parametrically, and Stoffer and Wall \cite{10.2307/2290521} developed a bootstrap algorithm for the maximum-likelihood estimator of state space models. Rodriguez and Ruiz tried to find the prediction interval for future observations through bootstrap \cite{https://doi.org/10.1111/j.1467-9892.2008.00604.x} but did not provide a theoretical justification on why their algorithm worked.

We are interested in the aforementioned problems, but would like to set constraints on the state space models and explain why we need the constraints at first. Example \ref{examp1} shows that two different state space models can generate the same observations. Specifically, a state space model can not be identified from observations if 1. distributions of state noise and measurement noise are unknown ((S1) and (S2)), 2. state transition matrix and measurement transition matrix are unknown ((S1) and (S3)) and 3. the dimension of state vectors are larger than the dimension of observations((S1) and (S4)). After excluding these situations, we may find methods to consistently estimate and predict a state space model.

\begin{example}
Consider the following 4 state space models
\begin{equation}
\begin{aligned}
\text{(S1)}
\begin{cases}
X_{n + 1} = \epsilon_{n + 1},\ \epsilon_1\sim\mathcal{N}(0, 1)\\
Y_n = X_n + \eta_n,\ \eta_1\sim\mathcal{N}(0, 2)
\end{cases}
\text{(S2)}
\begin{cases}
X_{n + 1} = \epsilon_{n+1},\ \epsilon_1 \sim\mathcal{N}(0, 2)\\
Y_n = X_n + \eta_n,\ \eta_1\sim\mathcal{N}(0, 1)
\end{cases}\\
\text{(S3)}
\begin{cases}
X_{n + 1} = \epsilon_{n+1},\ \epsilon_1 \sim\mathcal{N}(0, 4)\\
Y_n = \frac{1}{2} X_n + \eta_n,\ \eta_1\sim\mathcal{N}(0, 2)
\end{cases}
\text{(S4)}
\begin{cases}
X_{n+1} = \epsilon_{n+1},\ \epsilon_1\sim\mathcal{N}(0, I_2)\\
Y_n = [1, 0]X_{n+1} + \eta_n,\ \eta_1 \sim\mathcal{N}(0, 2)
\end{cases}
\end{aligned}
\end{equation}
Here $\mathcal{N}(0, a)$ denotes multivariate normal distribution with mean $0$ and covariance matrix $a$, and if $a$ is a scalar, then it denotes normal distribution with mean $0$ and variance $a$. $I_2$ is the two dimensional identity matrix. $\epsilon_i$ and $\eta_i$, $i = ...,-1,0,1,...$ are independent and identically distributed and $\epsilon_i$ is independent of $\eta_j$ for any $i,j$. These four models are different but their observations $Y_i,\ i = 1,2,...,n$ have the same joint distributions. In other words, a set of observations $Y_i,\ i = 1,2,...,n$ can be generated from (S1) to (S4).
\label{examp1}
\end{example}

In this paper, we assume that the measurement transition matrix and the distribution of measurement noise are known in addition to a set of observations from a state space model. Based on these information, we will provide a method to consistently estimate the state transition matrix and the distribution of state noise. Moreover, focus on filtering and prediction problems, we will provide consistent prediction intervals for unobserved state vectors and future observations.

The rest of this paper is organized as follows. In chapter \ref{Prelim} we introduce the problem setting and the frequently used notations and assumptions. In chapter \ref{chpEst} we present consistent estimators of the state transition matrix and the distribution of state noise and the related asymptotic results. In chapter \ref{cp5} we introduce algorithms to construct prediction intervals and explain why these algorithms maintain consistency. We perform finite sample simulations in chapter \ref{cp6} and make conclusions in chapter \ref{Conc}. Technique details are presented in the appendix.
\section{Preliminary}
\label{Prelim}
In the following of this paper, we focus on the linear state space model
\begin{equation}
\begin{aligned}
X_{n + 1} = AX_n + \epsilon_{n + 1}\\
Y_n = BX_n + \eta_n,\ \ n = ...,-1,0,1,...
\end{aligned}
\label{MainModel}
\end{equation}
According to chapter 1.3.3 in \cite{HMM}, we call $A$ the state transition matrix, $B$ the measurement transition matrix, $\epsilon_n$ the state noise and $\eta_n$ the measurement noise, $n = ...,-1,0,1,...$. $A$ and $B$ are $d\times d$ matrix and $\eta_n, \epsilon_n, X_n$ and $Y_n$ are $d$-dimensional random vectors. We respectively denote $f_\epsilon,\ f_\eta$ and $f_X$ as the density of $\epsilon_1,\ \eta_1$ and $X_1$. for a function $f:\mathbf{R}^d\to\mathbf{C}$, we define its fourier transformation and inverse Fourier transformation as
\begin{equation}
\begin{aligned}
\mathcal{F}f(x) = \int_{\mathbf{R}^d}f(y)\exp(iy^Tx)dy\\
\mathcal{F}^{-1}f(x) = \frac{1}{(2\pi)^d}\int_{\mathbf{R}^d}f(y)\exp(-iy^Tx)dy
\end{aligned}
\end{equation}
In addition, for a vector $x = (x_1,...,x_d)^T\in\mathbf{R}^d$, we denote its Euclidean norm as $\Vert x\Vert_2 = \sqrt{\sum_{i = 1}^d x^2_i}$, its infinity norm as $\Vert x\Vert_\infty = \max_{i = 1,2,...,d}\vert x_i\vert$. In addition, for a given non-singular matrix $M$, we denote the vector norm $\Vert x\Vert_M = \Vert Mx\Vert_\infty$, from the equivalence of norms(e.g. Exercise 5.6 in \cite{Folland:1706460}), we can find constants $c_M, C_M > 0$ such that $c_M\Vert x\Vert_M\leq\Vert x\Vert_\infty\leq C_M\Vert x\Vert_M$ for any $x\in\mathbf{R}^d$. For a matrix $M$, we define its $2$-norm as $\Vert M\Vert_2 = \sup_{x\in\mathbf{R}^d, \Vert x\Vert_2 = 1}\Vert Mx\Vert_2$. For a function $f:\mathbf{R}^d\to\mathbf{C}$, we define its $L_2$ norm as $\Vert f\Vert_{L_2} = \sqrt{\int_{\mathbf{R}^d}\vert f(x)\vert^2 dx}$. Similar as chapter 1.5.1 in \cite{mathematicsStat}, for two sequences of real numbers $a_n$ and $b_n > 0$, $n = 1,2,...$, we say that $a_n = O(b_n)$ if there exists a constant $C>0$ such that $a_n \leq C b_n$ for $\forall n$ and $a_n = o(b_n)$ if $a_n/b_n\to 0$ as $n\to\infty$. For two sequences of random variables $X_n$ and $Y_n > 0$, $n = 1,2,...$, we say that $X_n = O_p(Y_n)$ if for any given $\delta>0$, there exists a constant $C_\delta > 0$ such that $\sup_{n}Prob\left(\vert X_n\vert\geq C_\delta Y_n\right)\leq \delta$ and $X_n = o_p(Y_n)$ if $X_n/Y_n\to_p 0$ as $n\to\infty$. We denote $\to_p$ as convergence in probability. Other symbols will be defined before being used. With these notations, we introduce the main assumptions of this paper.

\textbf{Assumptions}

1). For $\tau = \epsilon, \eta$, in model \eqref{MainModel} $\tau_n\in\mathbf{R}^d, n = ...,-1,0,1,...$ are independent and identically distributed, and the distribution of $\tau_1$ is absolutely continuous with respect to Lebesgue measure with density $f_\tau$. $\mathbf{E}\vert\tau_{1,j}\vert < \infty$ for $j = 1,2,...,d$ and $\mathbf{E}\tau = 0$. $\epsilon_n, n =...,-1,0,1,...$ are independent of $\eta_m, m = ...,-1,0,1,...$. $d\times d $ matrix $A$, $B$ and the covariance matrix $\Sigma$ of $\epsilon_1$ are non-singular. In addition, we assume that $\Vert A\Vert_2 < 1$, the distribution of $X_n$ is stationary(like definition 1.3.3 in \cite{time_series_analysis}). We also assume that the measurement transition matrix $B$ and the density of measurement noise $f_\eta$ are known in addition to observations $Y_i,\ i = 1,2,...,n$.

2). For $\tau = \epsilon, \eta$, we assume that
\begin{equation}
\begin{aligned}
\sum_{j = 1}^d\int_{\mathbf{R}^d}x_j^4f_\tau(x_1,...,x_d)dx_1...dx_d < \infty\\
\int_{\mathbf{R}^d}\vert\mathcal{F}f_\tau(x_1,...x_d)\vert dx_1...dx_d < \infty
\end{aligned}
\end{equation}

3). $G:\mathbf{R}\to\mathbf{R}$ is a function such that $\int_{\mathbf{R}}\vert G(x)\vert dx <\infty$, and $\mathcal{F}G$'s support belongs to $[-a,a]$ for a constant $a > 1$, $0\leq \mathcal{F}G \leq 1$ and $\mathcal{F}G = 1$ on $[-1, 1]$. For each bandwidth $h > 0$, we define the kernel function in $\mathbf{R}^d$ as
\begin{equation}
K_h(x_1,...,x_d) = \prod_{i = 1}^d \frac{1}{h}G\left(\frac{x_i}{h}\right)\Rightarrow \mathcal{F}K_h(x_1,...,x_d) = \prod_{i = 1}^d\mathcal{F}G(hx_i)
\end{equation}

4). One of the two following conditions happens:

4.a) (Ordinary smooth) There exists constants $b,\beta,c > 0$ such that
\begin{equation}
\vert\mathcal{F}f_\eta(x_1,...x_d)\vert\geq c\prod_{j = 1}^d(x_j^2 + 1)^{-\beta / 2},\ \sum_{j = 1}^d\int_{\mathbf{R}^d}\vert\mathcal{F}f_\epsilon(x_1,...,x_d)\vert^2(1 + x_j^2)^{b / 2}dx_1...dx_d < \infty
\end{equation}
and the bandwidth $h = h(n)$ satisfies $h = o(1)$ and $\frac{1}{nh^{4d\beta + d + 2}} = o(1)$.

4.b) (Super smooth) There exists constants $c,\beta,\gamma, b, r > 0 $ such that
\begin{equation}
\vert\mathcal{F}f_\eta(x_1,...,x_d)\vert\geq c\prod_{j = 1}^d\exp(-\gamma\vert x_j\vert^\beta),\ \sum_{j = 1}^d \int_{\mathbf{R}^d}\vert\mathcal{F}f_\epsilon(x_1,...,x_d)\vert^2\exp(r\vert x_j\vert^b)dx_1...dx_d < \infty
\end{equation}
and the bandwidth $h = h(n)$ satisfies $h = o(1)$ and $\frac{1}{h\log^{1/\beta}(n)} = o(1)$.

Classifying distributions through their smoothness is common in the deconvolution literatures such as Fan \cite{https://doi.org/10.2307/3315465}, Fan \cite{10.2307/24304026} and Cao \cite{DeconvolutionCumu}. Assumption 4) is similar with conditions used by Comte and Lacour \cite{AIHPB_2013__49_2_569_0} for solving multivariate deconvolution problems. In remark \ref{remark1}, we derive some important corollaries from aforementioned assumptions.
\begin{remark}
\label{remark1}
By letting $\tau = \epsilon, \tau$, from theorem 8.22 and 8.26 in \cite{Folland:1706460} and assumption 2), we know that $\mathcal{F}f_\tau$ is continuous and vanishes at infinity, and $\mathcal{F}\mathcal{F}^{-1}f_\tau = \mathcal{F}^{-1}\mathcal{F}f_\tau = f_\tau$ almost surely. Thus, without loss of generality, we may assume that $f_\tau$ is continuous(which implies uniformly continuous). According to Plancherel's theorem, $\Vert\mathcal{F}f_\epsilon\Vert_{L_2} = (2\pi)^{d/2}\Vert f_\epsilon\Vert_{L_2}$.

The second thing we would like to emphasize is that $G$ defined in assumption 3) exists. For example, we can set $a = 2$ and define function $G$ as
\begin{equation}
G(x)=\frac{1}{\pi x^2}\left(\cos(x)-\cos(2x)\right)\Rightarrow \mathcal{F}G(x)=
\begin{cases}
x+2\ \ x\in[-2,-1]\\
1\ \ x\in[-1,1]\\
2-x\ \ x\in[1,2]\\
0\ \ \text{Otherwise}
\end{cases}
\label{Kernel}
\end{equation}

For any Borel-measurable set $S$ with Lebesgue measure $0$, define set $S - Ax = \{y - Ax|y\in S\},\ \forall x\in\mathbf{R}^d$, $S - Ax$ has Lebesgue measure $0$ as well, so $Prob\left( X_1\in S\right) = \mathbf{E}Prob\left(\epsilon_1\in S - AX_0|X_0\right) = 0$, which implies that the distribution of $X_1$ is absolutely continuous with respect to Lebesgue measure and the density $f_X$ is well-defined. From assumption 1), we have for any $x\in\mathbf{R}^d$,
\begin{equation}
\mathcal{F}f_X(x) = \mathcal{F}f_X(A^Tx)\mathcal{F}f_\epsilon(x)\Rightarrow\vert \mathcal{F}f_X(x)\vert\leq \vert\mathcal{F}f_\epsilon(x)\vert
\end{equation}
Therefore, assumption 2) and 4) are also valid for $\mathcal{F}f_X$.
\end{remark}
\section{Semi-parametric estimation in state space models}
\label{chpEst}
In this chapter we will find estimators for state transition matrix $A$ and the density of state noises $f_\epsilon$ as well as the density of state vectors $f_X$. From \eqref{MainModel}, by defining $\mathcal{F}f_Y(x) = \mathbf{E}\exp(ix^TY_1)$,
\begin{equation}
\begin{aligned}
\mathcal{F}f_Y(x) = \mathcal{F}f_X(B^Tx)\mathcal{F}f_\eta(x),\ \
Y_{n+1} - BAB^{-1}Y_n = \eta_{n+1} + B\epsilon_{n+1} - BAB^{-1}\eta_n\\
\Rightarrow \mathcal{F}f_{Y_2 - BAB^{-1}Y_1}(x) = \mathcal{F}f_\eta(x)\mathcal{F}f_\epsilon(B^Tx)\mathcal{F}f_\eta(-B^{-1T}A^TB^Tx)
\end{aligned}
\end{equation}
Since $\mathcal{F}f_Y$ can be consistently estimated from observations, we may apply deconvolution to achieve this goal. We refer Fan \cite{10.2307/24304026}, Fan \cite{https://doi.org/10.2307/3315465} and Devroye \cite{10.2307/3314852} as an introduction. Masry \cite{MASRY199353} considered deconvolution problem for dependent data and Comte and Lacour \cite{AIHPB_2013__49_2_569_0} tried to perform deconvolution for random vectors. We present the estimators and summarize the asymptotic results in theorem \ref{thm1}.
\begin{theorem}
Suppose assumptions 1) to 4),

1. By defining the estimator $\widehat{A}$ as
\begin{equation}
\widehat{A}=(\sum_{k=3}^n B^{-1}Y_kY_{k-2}^TB^{-1T})(\sum_{k=3}^n B^{-1}Y_{k-1}Y_{k-2}^TB^{-1T})^+
\label{Esti_A}
\end{equation}
$\Vert\widehat{A} - A\Vert_2 = O_p(1/\sqrt{n})$, here we denote $A^+$ as the pseudo-inverse of matrix $A$.

2. By defining the estimator $\widehat{f}_X(x),\ x\in\mathbf{R}^d$ as
\begin{equation}
\widetilde{f}_X(x) = \frac{1}{n\times(2\pi)^d}\sum_{j = 1}^n \int_{\mathbf{R}^d}\frac{\mathcal{F}K_h(y)}{\mathcal{F}f_\eta(B^{-1T}y)}\exp(iy^TB^{-1}Y_j - iy^Tx)dy,\ \widehat{f}_X(x) = \max(\Re(\widetilde{f}_X(x)), 0)
\end{equation}
with $\Re(x)$ denotes the real part of a complex number $x$, we have $\mathbf{E}\Vert\widehat{f}_X - f_X\Vert_{L_2}^2 = o(1)$.

3. By defining the estimator $\widehat{f}_\epsilon(x),\ \forall x\in\mathbf{R}^d$ as
\begin{equation}
\begin{aligned}
\widetilde{f}_\epsilon(x) = \frac{1}{(n - 1)\times (2\pi)^d}\sum_{j = 1}^{n - 1}\int_{\mathbf{R}^d}\frac{\mathcal{F}K_h(y)}{\mathcal{F}f_\eta(B^{-1T}y)\mathcal{F}f_\eta(-B^{-1T}\widehat{A}^Ty)}\exp(iy^TB^{-1}Y_{j+1} - iy^T\widehat{A}B^{-1}Y_j - iy^Tx)dy\\
\widehat{f}_\epsilon(x) = \max(\Re(\widetilde{f}_\epsilon(x)), 0)
\end{aligned}
\label{feps}
\end{equation}
we have $\Vert\widehat{f}_\epsilon - f_\epsilon\Vert_{L_2} = o_p(1)$
\label{thm1}
\end{theorem}

\begin{remark}
Even though we present estimators $\widetilde{f}_X$ and $\widetilde{f}_\epsilon$ through multiple integral, in practice we can derive them through Monte Carlo integration. Notice that $\mathcal{F}K_h(y) = 0$ for $\Vert y\Vert_\infty > a / h$, in practice for a given positive integer $R$, we can generate $R\times d$ independent and identically distributed random variables $\zeta_{jk}, j = 1,2,...,R, k = 1,2,...,d$ whose marginal distribution is uniform distribution on $[-a/h, a/h]$. Then we define $\zeta_i = (\zeta_{i1},...,\zeta_{id})^T$ and derive
\begin{equation}
\begin{aligned}
I(\widetilde{f}_X(x)) = \left(\frac{a}{h\pi}\right)^d\frac{1}{nR}\sum_{k = 1}^R\sum_{j = 1}^n\frac{\mathcal{F}K_h(\zeta_k)}{\mathcal{F}f_\eta(B^{-1T}\zeta_k)}\exp(i\zeta_k^TB^{-1}Y_j - i\zeta_k^Tx)\\
I(\widetilde{f}_\epsilon(x)) = \left(\frac{a}{h\pi}\right)^d\frac{1}{(n - 1)R}\sum_{k = 1}^R\sum_{j = 1}^{n - 1}\frac{\mathcal{F}K_h(\zeta_k)\exp(i\zeta_k^TB^{-1}Y_{j+1} - i\zeta_k^T\widehat{A}B^{-1}Y_j - i\zeta_k^Tx)}{\mathcal{F}f_\eta(B^{-1T}\zeta_k)\mathcal{F}f_\eta(-B^{-1T}\widehat{A}^T\zeta_k)}
\end{aligned}
\label{HandCal}
\end{equation}
from the strong law of large number, when $\widehat{A}$ is non-singular, from the equivalence of norm, we have $\lim_{R\to\infty}I(\widetilde{f}_X(x)) = \mathbf{E}I(\widetilde{f}_X(x))|Y_1,...,Y_n = \widetilde{f}_X(x)$ and $\lim_{R\to\infty}I(\widetilde{f}_\epsilon(x)) = \mathbf{E}I(\widetilde{f}_\epsilon(x)) |Y_1,...,Y_n = \widetilde{f}_\epsilon(x)$ almost surely. Therefore, by selecting sufficiently large $R$, we may use \eqref{HandCal}, which is easily calculated, to approximate $\widetilde{f}_X$ and $\widetilde{f}_\epsilon$.
\end{remark}

As a special case, if $A = 0$ and $B = I_d$, the $d$ dimensional identity matrix, then the state space model \eqref{MainModel} degenerates to the usual de-convolution model $Y_n = \epsilon_{n} + \eta_n$, and the estimator $\widetilde{f}_X$ coincides with the frequently used estimator in de-convolution setting.

We emphasize an important corollary from theorem \ref{thm1}. Suppose a function $g$ satisfies $\Vert g\Vert_{L_2} < \infty$, by defining the convolution $f\star g(x) = \int_{\mathbf{R}^d}f(y)g(x - y)dy$, for $\tau = X, \epsilon$, from Cauchy's inequality
\begin{equation}
\sup_{x\in\mathbf{R}^d}\vert\widehat{f}_\tau(x) - f_\tau(x)\vert\leq \sup_{x\in\mathbf{R}^d}\int_{\mathbf{R}^d}\vert\widehat{f}_\tau(y) - f_\tau(y)\vert\times \vert g(x - y)\vert dy\leq \Vert g\Vert_{L_2}\times \Vert\widehat{f}_\tau - f_\tau\Vert_{L_2} = o_p(1)
\end{equation}
\section{Prediction intervals in state space models}
\label{cp5}
Prediction is not a new topic in statistical inference but has seen resurgence with the advent of statistical learning (Hastie et.al.\cite{ESl}). We refer Geisser \cite{PrdInference} for a comprehensive introduction and Politis \cite{model_free} for a more recent exposition. Apart from a precise predictor, statisticians always use prediction intervals to quantify uncertainty in a prediction \cite{NEURIPS2019_5103c358}. In this chapter, we will provide an algorithm to find the prediction intervals for state vectors and future observations and theoretically justify why the presented algorithm generates consistent prediction intervals. We adopt definition 2.4.1 in \cite{model_free} to define a consistent prediction interval.

\begin{definition}[Consistent prediction interval]
Suppose there are $n$ observations $Y_1,...,Y_n$ and the set $\Gamma = \Gamma(Y_1,...,Y_n)$ as well as the new random variable $Z_n$ satisfy
\begin{equation}
Prob\left(Z_n\in\Gamma\right)\to 1-\alpha\ \ \text{as }n\to\infty
\end{equation}
then we call the set $\Gamma$ the $1-\alpha$ consistent prediction interval for the random variable $Z_n$. Here $0< \alpha < 1$.
\label{def1}
\end{definition}
In this paper, we are interested in $Z_n = X_n$(filtering) and  $Z_n = X_{n+1},\ Y_{n+1}$(prediction).

\begin{example}[Consistent prediction intervals for standard state space model] We adopt the definition of a standard state space model in chapter 1.3.3, \cite{HMM}, i.e. $A,\ B$ and the covariance matrix $\Sigma$ of $\epsilon_1$ as well as the covariance matrix $\Pi$ of $\eta_1$ are known. $\epsilon_1,\ \eta_1$ are multivariate normal distributed. In addition, we assume $X_1 = [0, 0,...,0]^T$ and $Y_0 = [0, 0,...,0]^T$. From Kalman prediction (proposition 12.2.2 in \cite{time_series_analysis}), by defining $\widehat{X}_n$ as the Kalman one-step predictor for $X_n$, we may derive the error covariance matrix $\Omega_n = \mathbf{E}(X_n - \widehat{X}_n)(X_n - \widehat{X}_n)^T$. Since $\mathbf{E}X_n = 0$ and $\mathbf{E}\widehat{X}_n = 0$ (we can prove this result by $\mathbf{E}\widehat{X}_1 = 0$ and induction), if $\Omega_n$ is non-singular, one of the $1-\alpha$ consistent prediction interval for $X_n$ can be given by $(X_n - \widehat{X}_n)^T\Omega_n^{-1}(X_n - \widehat{X}_n)\leq c_{\chi^2_d, 1-\alpha}$, here $c_{\chi^2_d, 1-\alpha}$ is the $1-\alpha$ quantile of $\chi^2$ distribution with degree of freedom $d$. Since $X_{n+1} - A\widehat{X}_n = A(X_n - \widehat{X}_n) + \epsilon_{n+1}$, the $1-\alpha$ consistent prediction interval for $X_{n+1}$ can be given by $(X_{n+1} - A\widehat{X}_n)^T(A\Omega_n A^T + \Sigma)^{-1}(X_{n+1} - A\widehat{X}_n)\leq c_{\chi^2_d, 1-\alpha}$. Since $Y_{n + 1} - BA\widehat{X}_n = \eta_{n+1} + B\epsilon_{n+1} + BA(X_n - \widehat{X}_n)$, the $1-\alpha$ consistent prediction interval for $Y_{n+1}$ can be derived by
\begin{equation}
(Y_{n+1} - BA\widehat{X}_n)^T(\Pi + B\Sigma B^T + BA\Omega_n A^TB^T)^{-1}(Y_{n+1} - BA\widehat{X}_n)\leq c_{\chi^2_d, 1-\alpha}
\end{equation}
Consistency of the aforementioned prediction intervals relies on two facts: 1. the distribution of $X_n - \widehat{X}_n$ only depends on the covariance matrix $\Omega_n$, and 2. linear combinations of random variables preserve normality. These two facts fail if $\epsilon_1$ and $\eta_1$ are not normally distributed.
\label{EXPZ}
\end{example}

Kalman recursion (chapter 12.2 in \cite{time_series_analysis}) is one of the most popular methods to deal with state space models. However, when the state noise and the measurement noise are not normal, distribution of the predictor will be hard to derive. Capp\'{e} et.al. \cite{HMM} introduced how to filter and smooth a general state space model but they assumed nothing needed to be estimated. We use another method to derive predictors. If $\eta_i = 0$, then model \ref{MainModel} degenerates to a vector autoregressive model.  In this situation, according to Kim \cite{KIM1999393}, predictors for $X_n$, $X_{n+1}$ and $Y_{n+1}$ can respectively be $B^{-1}Y_n$, $\widehat{A}B^{-1}Y_n$ and $B\widehat{A}B^{-1}Y_n$. When $\eta_i \neq 0$, we still can use these predictors but the distributions of predictive roots(definition 2.4.2 in \cite{model_free}) are different from vector autogressive models.

\begin{lemma}
\label{lemma2}
Suppose $Y_i, i =1,2,...,n$ are generated from state space model \ref{MainModel} and assumption 1) is satisfied, then for $\forall 0 < x\in\mathbf{R}$

1. Filtering:
\begin{equation}
Prob\left(\Vert X_n - B^{-1}Y_n\Vert_\infty\leq x\right) = \int_{\Vert B^{-1}y\Vert_\infty\leq x}f_\eta(y)dy
\label{Z1}
\end{equation}

2. Predicting the state vector:
\begin{equation}
Prob\left(\Vert X_{n+1} - AB^{-1}Y_n\Vert_\infty\leq x\right) = \int_{\mathbf{R}^d}\left(\int_{\Vert z\Vert_\infty\leq x}f_\epsilon(z + AB^{-1}y)dz\right)f_\eta(y)dy
\end{equation}

3. Predicting the future observation:
\begin{equation}
Prob\left(\Vert Y_{n+1} - BAB^{-1}Y_n\Vert_\infty\leq x\right) = \frac{1}{\vert\det(B)\vert}\int_{\mathbf{R}^d\times\mathbf{R}^d}\left(\int_{\Vert y\Vert_\infty\leq x}f_\epsilon(B^{-1}y + AB^{-1}w - B^{-1}z)dy\right)f_\eta(w)f_\eta(z)dwdz
\label{Z3}
\end{equation}
\end{lemma}

Algorithm \ref{alg1} performs binary searches on the estimated cumulative distribution functions of predictive roots $X_n - B^{-1}Y_n$, $X_{n+1} - \widehat{A}B^{-1}Y_n$ and $Y_{n+1} - B\widehat{A}B^{-1}Y_n$ to find $1-\alpha$ quantiles. Validness of algorithm \ref{alg1} can be assured if the estimated distributions converge to the distributions of predictive roots. We prove this result in theorem \ref{MainThm} and corollary \ref{MainCoro}.
\IncMargin{1em}

\begin{algorithm}[htbp]
\DontPrintSemicolon
\SetKwData{Left}{left}\SetKwData{This}{this}\SetKwData{Up}{up}
\SetKwFunction{Union}{Union}\SetKwFunction{FindCompress}{FindCompress}
\SetKwInOut{Input}{input}\SetKwInOut{Output}{output}
\Input{Observations $Y_1,...,Y_n\in\mathbf{R}^d$ from model \eqref{MainModel}, kernel function $G$, bandwidth $h$, matrix $B$, density $f_\eta$, nominal coverage probability $1-\alpha\in(0,1)$, number of replicates $R$, tolerance $\varepsilon$}

\Output{Consistent prediction intervals for $X_{n}$, $X_{n+1}$ and $Y_{n+1}$}
\BlankLine

Derive statistics $\widehat{A}$, $\widehat{f}_\epsilon$ defined in theorem \ref{thm1}\;

\For{operation $\mathcal{M}\in \{\mathcal{H},\mathcal{N},\mathcal{G}\}$}{
Set $x_{\mathcal{M}, l} = 0,\ x_{\mathcal{M}, h} = 1$\;
\lWhile{$\mathcal{M}(x_{\mathcal{M}, h}) < 1-\alpha$}{
$x_{\mathcal{M}, h} = 2\times x_{\mathcal{M}, h}$\;
}
\While{$x_{\mathcal{M}, h} - x_{\mathcal{M}, l} > \varepsilon$}{
Set $x_{\mathcal{M}, m} = \frac{x_{\mathcal{M}, h} + x_{\mathcal{M}, l}}{2}$, calculate $s = \mathcal{M}(x_{\mathcal{M}, m})$\;
\lIf{$s < 1-\alpha$}{$x_{\mathcal{M}, l} = x_{\mathcal{M}, m}$}
\lElse{$x_{\mathcal{M}, h} = x_{\mathcal{M}, m}$}
}
Set $x_{\mathcal{M}} = \frac{x_{\mathcal{M}, h} + x_{\mathcal{M}, l}}{2}$
}

The prediction intervals for $X_n,\ X_{n+1}$ and $Y_{n+1}$ respectively is given by $\Vert X_n - B^{-1}Y_n\Vert_\infty\leq x_{\mathcal{H}}$, $\Vert X_{n+1} - \widehat{A}B^{-1}Y_n\Vert_\infty\leq x_{\mathcal{N}}$
and $\Vert Y_{n+1} - B\widehat{A}B^{-1}Y_n\Vert_\infty\leq x_{\mathcal{G}}$
\BlankLine
\BlankLine
\textbf{Operation} $\mathcal{H}(x)$:\;
Generate $R$ independent and identically distributed random vectors $\zeta_j,\ j = 1,2,...,R$ with marginal density $f_\eta$\;
Derive the sample mean $\tau = \frac{1}{R}\sum_{j = 1}^R\mathbf{1}_{\Vert B^{-1}\zeta_j\Vert_\infty\leq x}$\;
\Output{$\tau$}
\BlankLine
\BlankLine
\textbf{Operation} $\mathcal{N}(x)$:\;
Generate $R$ independent and identically distributed random vectors $\zeta_j,\ j = 1,2,..., R$ with marginal density $f_\eta$ and $R\times d$ independent and identically distributed random variables $\kappa_{jk},$ $j = 1,2,...,R$ and $k = 1,2,...,d$ with uniform marginal distribution on $[-x, x]$, denote $\kappa_j = [\kappa_{j1},...,\kappa_{jd}]^T$ for $j = 1,2,...,R$\;
Derive the sample mean $\tau = \frac{1}{R}\sum_{j = 1}^R\widehat{f}_\epsilon(\kappa_j + \widehat{A}B^{-1}\zeta_j)$\;
\Output{$(2x)^d\times\tau$}
\BlankLine
\BlankLine
\textbf{Operation} $\mathcal{G}(x)$:\;
Generate $2R$ independent and identically distributed random vectors $\zeta_{s,j}$, $s = 1,2$ and $j = 1,2,...,R$ with marginal density $f_\eta$, and generate $R\times d$ independent and identically distributed random variables $\kappa_{j,k}$, $j = 1,...,R$ and $k = 1,...,d$, denote $\kappa_{j} = [\kappa_{j,1},...,\kappa_{j,d}]^T$\;
Derive $\tau = \frac{1}{R}\sum_{j = 1}^R\widehat{f}_\epsilon(B^{-1}\kappa_j + \widehat{A}B^{-1}\zeta_{1,j} - B^{-1}\zeta_{2,j})$\;
\Output{$\frac{(2x)^d}{\vert\det(B)\vert}\times\tau$}
\label{alg1}
\caption{Constructing prediction intervals for $X_n$, $X_{n+1}$ and $Y_{n+1}$}
\end{algorithm}
\DecMargin{1em}

\begin{remark}
From the strong law of large number, by defining $\underline{Y}_n = (Y_1,...,Y_n)$, almost surely for any given $x\in\mathbf{R}^d$,
\begin{equation}
\begin{aligned}
H(x) = \lim_{R\to\infty}\frac{1}{R}\sum_{j = 1}^R\mathbf{1}_{\Vert B^{-1}\zeta_j\Vert_\infty\leq x} = Prob\left(\Vert B^{-1}\zeta_1\Vert_\infty\leq x|\underline{Y}_n\right) = \int_{\Vert B^{-1}y\Vert_\infty\leq x}f_\eta(y)dy\\
N(x) = (2x)^d\lim_{R\to\infty}\frac{1}{R}\sum_{j = 1}^R\widehat{f}_\epsilon(\kappa_j + \widehat{A}B^{-1}\zeta_j) = (2x)^d\mathbf{E}\widehat{f}_\epsilon(\kappa_1 + \widehat{A}B^{-1}\zeta_1)|\underline{Y}_n\\
= \int_{\Vert z\Vert_\infty\leq x, y\in\mathbf{R}^d}\widehat{f}_\epsilon(z + \widehat{A}B^{-1}y)f_\eta(y)dzdy\\
G(x) = \frac{(2x)^d}{\vert\det(B)\vert}\lim_{R\to\infty}\sum_{j = 1}^R\widehat{f}_\epsilon(B^{-1}\kappa_j + \widehat{A}B^{-1}\zeta_{1,j} - B^{-1}\zeta_{2,j})
= \frac{(2x)^d}{\vert\det(B)\vert}\mathbf{E}\widehat{f}_\epsilon(B^{-1}\kappa_1 + \widehat{A}B^{-1}\zeta_{1,1} - B^{-1}\zeta_{2,1})|\underline{Y}_n\\
= \frac{1}{\vert\det(B)\vert}\int_{\mathbf{R}^d\times\mathbf{R}^d}\left(\int_{\Vert y\Vert_\infty\leq x}\widehat{f}_\epsilon(B^{-1}y + \widehat{A}B^{-1}w - B^{-1}z)dy\right)f_\eta(w)f_\eta(z)dwdz
\end{aligned}
\end{equation}
which are the desired estimators.

Compare to classical de-convolution literatures like \cite{https://doi.org/10.2307/3315465} or \cite{AIHPB_2013__49_2_569_0}, algorithm \ref{alg1} avoids calculating multidimensional integrations, which is hard to manipulate. There are several methods to generate random vectors with desired density $f_\eta$, e.g., Metropolis-Hastings algorithm \cite{MCMC} and Accept-Reject algorithm \cite{HMM}.

From \eqref{feps}, we know that $\widehat{f}_\epsilon$ is bounded, so from dominated convergence theorem $H, N$ and $G$ are continuous and increasing. By defining $M = H, N,G$, if $M(0) < 1 - \alpha$ and $\lim_{x\to\infty}M(x) > 1-\alpha$, from intermediate value theorem there exists $x^*_{\mathcal{M}}$ such that $M(x^*_{\mathcal{M}}) = 1 - \alpha$. Since $x^*_{\mathcal{M}}$ can be approximated in algorithm \ref{alg1} through setting $R\to\infty$ and $\varepsilon\to 0$, we only consider $M$ and $x^*_{\mathcal{M}}$ in theoretical justifications.
\label{remark2}
\end{remark}

\begin{theorem}
\label{MainThm}
Suppose assumptions 1) to 4), then for arbitrary given $x_0 > 0$
\begin{equation}
\begin{aligned}
\sup_{0 < x\leq x_0}\vert N(x) - Prob\left(\Vert X_{2} - AB^{-1}Y_1\Vert_\infty\leq x\right)\vert = o_p(1)\\
\sup_{0<x\leq x_0}\vert G(x) - Prob\left(\Vert Y_{2} - BAB^{-1}Y_1\Vert_\infty\leq x\right)\vert = o_p(1)
\end{aligned}
\label{ImportEq}
\end{equation}
Here $H,\ N$ and $G$ coincide with remark \ref{remark2}.
\end{theorem}

For a given $0 < \alpha < 1$, by choosing sufficiently large $x_0$ such that
\begin{equation}
Prob\left(\Vert X_{2} - AB^{-1}Y_1\Vert_\infty\leq x_0\right) > 1 - \alpha / 4\ \text{and}\  Prob\left(\Vert Y_{2} - BAB^{-1}Y_1\Vert_\infty\leq x_0\right) > 1-\alpha / 4
\label{Cond}
\end{equation}
and apply theorem \ref{MainThm} on $(0, x_0]$, the probability of $x_{\mathcal{N}}^*$ and $x_{\mathcal{G}}^*$ exist converges to $1$ as the sample size increases. In corollary \ref{MainCoro} we explain why $x_{\mathcal{M}}^*,\ \mathcal{M} = \mathcal{H}, \mathcal{N}, \mathcal{G}$, the asymptotic values of $x_\mathcal{M}$ in algorithm \ref{alg1}, contribute to the consistent prediction intervals.

\begin{corollary}
\label{MainCoro}
Suppose assumptions 1) to 4), for arbitrary given $0<\alpha<1$, by defining $x_{\mathcal{M}}^*,\mathcal{M} = \mathcal{H}, \mathcal{N}, \mathcal{G}$ as in remark \ref{remark2}
\begin{equation}
\begin{aligned}
Prob\left(\Vert X_n - B^{-1}Y_n\Vert_\infty\leq x_{\mathcal{H}}^*\right) - (1-\alpha) = 0\\
\vert Prob\left(\Vert X_{n+1} - \widehat{A}B^{-1}Y_n\Vert_\infty\leq x_{\mathcal{N}}^*\right) - (1 - \alpha)\vert = o(1)\\
\vert Prob\left(\Vert Y_{n+1} - B\widehat{A}B^{-1}Y_n\Vert_\infty\leq x^*_{\mathcal{G}}\right) - (1-\alpha)\vert = o(1)
\end{aligned}
\label{ResCoro}
\end{equation}
\end{corollary}
\section{Simulations and Numerical experiments}
\label{cp6}
In this chapter, we justify theoretical results in chapter \ref{chpEst} and \ref{cp5} through numerical experiments. We separate discussions into ordinary smooth and super smooth situations and design the following four models:

(O1) $X_n, Y_n$ satisfy $X_{n+1} = 0.8X_n + \epsilon_{n+1}$ and $Y_n = 1.0X_n + \eta_n$. $\epsilon_1 = \epsilon_{1,1} - \epsilon_{1,2}$, $\epsilon_{1,i},\ i = 1, 2$ have gamma distribution with shape $3 / 2$ and scale $1 / \sqrt{3}$, $\epsilon_{1,1}$ is independent of $\epsilon_{1,2}$. $\eta_1 = \eta_{1,1} - \eta_{1,2}$, $\eta_{1,i}.\ i = 1,2$ have gamma distribution with shape $1/2$ and scale $1.0$. $\eta_{1,1}$ is independent of $\eta_{1,2}$. In this situation, $\eta_1$ and $\epsilon_1$ have variance $1$ and $\beta = 1$. According to assumption 4.a), we choose the bandwidth $h = 1.0 / n^{1/8}$

(S1) We choose the same $A$ and $B$ as (O1) but $\epsilon_1$ and $\eta_1$ have normal distribution with mean $0$ and variance $1$. According to assumption 4.b), since $\beta = 2$, we choose $h = 1 / log^{0.1}(n)$.

(O2) $X_n, Y_n$ satisfy
\begin{equation}
\begin{aligned}
X_{n+1} =
\left[
\begin{matrix}
0.56 & -0.25 \\
0.25 & 0.45
\end{matrix}
\right]X_n + \left[
\begin{matrix}
0.979 & 0.204\\
0.204 & 0.979\\
\end{matrix}
\right]
\left[
\begin{matrix}
\epsilon_{n+1, 1}\\
\epsilon_{n+1, 2}
\end{matrix}
\right]\\
Y_n = \left[
\begin{matrix}
1.00 & -0.50\\
0.50 & 1.00\\
\end{matrix}
\right]X_n + \left[
\begin{matrix}
0.900 & 0.000\\
0.000 & 0.900
\end{matrix}
\right]
\left[
\begin{matrix}
\eta_{n, 1}\\
\eta_{n, 2}
\end{matrix}
\right]
\end{aligned}
\end{equation}
$\epsilon_{i,1}$ is independent of $\epsilon_{i,2}$, $\eta_{i, 1}$ is independent of $\eta_{i,2}$, and $\epsilon_{i,j}$, $\eta_{i,j},\ i = ...,-1,0,1,...,\ j = 1,2$ respectively has the same distribution as $\epsilon_1$ and $\eta_1$ in (O1). Since $\beta = 1$, we choose $h = 1.0 / n^{1/8}$.

(S2) We choose the same $A$ and $B$ as (O2) but $\epsilon_{i,j},\ j = 1,2$, $\eta_{i,j},\ j = 1,2$ have normal distribution with mean $0.0$ and variance $1.0$. Since $\beta = 2.0$, we choose $h = 1.0 / \log^{0.1}(n)$.

In the first experiment, we focus on evaluating the performance of presented estimators. Apart from $\widehat{A}, \widehat{f}_X$ and $\widehat{f}_\epsilon$, we also consider the estimator $\widehat{z}(x) = \int_{y\in\mathbf{R}^d}\widehat{f}_\epsilon(x + \widehat{A}B^{-1}y)f_\eta(y)dy$ on $z(x) = \int_{y\in\mathbf{R}^d}f_\epsilon(x + AB^{-1}y)f_\eta(y)dy$ and $ \widehat{g}(x) = \int_{\mathbf{R}^d\times\mathbf{R}^d}\widehat{f}_\epsilon(B^{-1}x +\widehat{A}B^{-1}w - B^{-1}z)f_\eta(w)f_\eta(z)dwdz$ on $g(x) = \int_{\mathbf{R}^d\times \mathbf{R}^d}f_\epsilon(B^{-1}x + AB^{-1}w - B^{-1}z)f_\eta(w)f_\eta(z)dwdz$. According to lemma \ref{lemma2}, if $\widehat{z}$ is close to $z$ and $\widehat{g}$ is close to $g$, algorithm \ref{alg1} should be able to find the consistent prediction interval for $X_{n + 1}$ and $Y_{n+1}$.

\begin{figure}[htbp]
\subfigure[$f_\epsilon$, model (O1), size 500]{
\includegraphics[width= 2.0in]{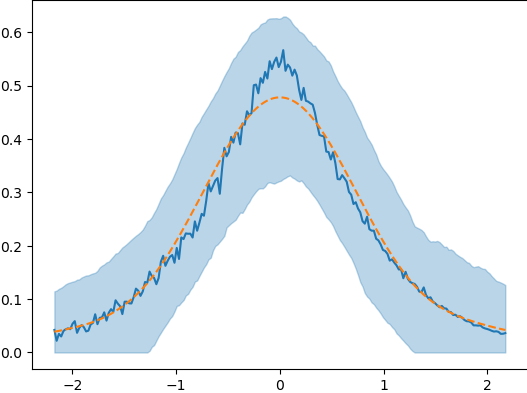}
\label{2case}
}
\subfigure[$z$, model (O1), size 500]{
\includegraphics[width = 2.0in]{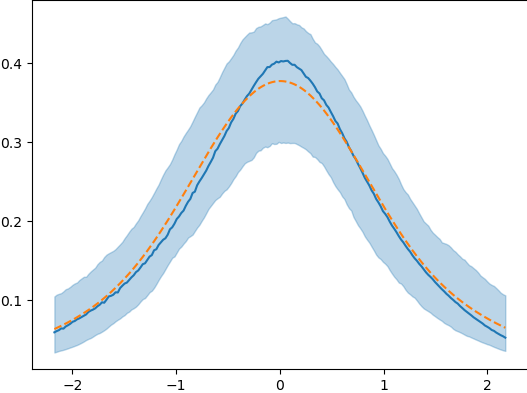}
\label{2case}
}
\subfigure[$g$, model (O1), size 500]{
\includegraphics[width = 2.0in]{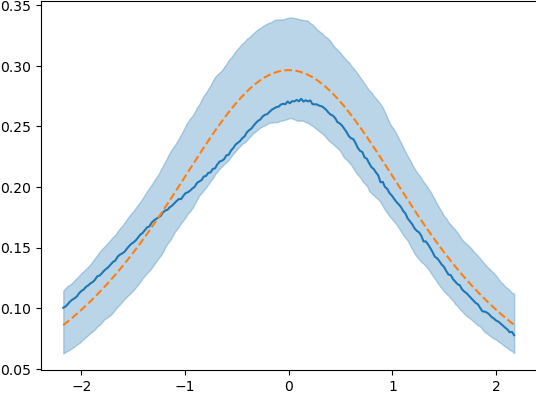}
\label{2case}
}

\subfigure[$f_\epsilon$, model (S1), size 500]{
\includegraphics[width= 2.0in]{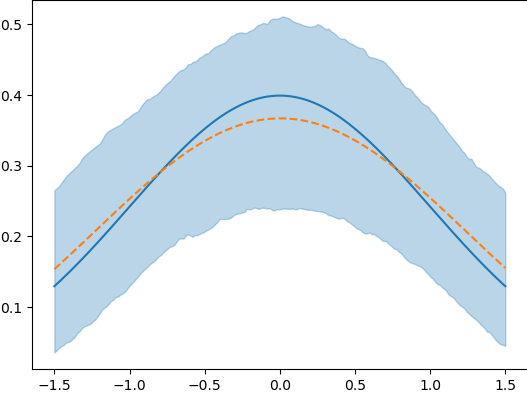}
\label{2case}
}
\subfigure[$z$, model (S1), size 500]{
\includegraphics[width = 2.0in]{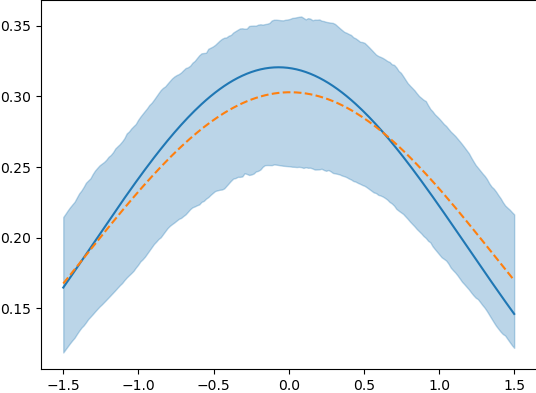}
\label{2case}
}
\subfigure[$g$, model (S1), size 500]{
\includegraphics[width = 2.0in]{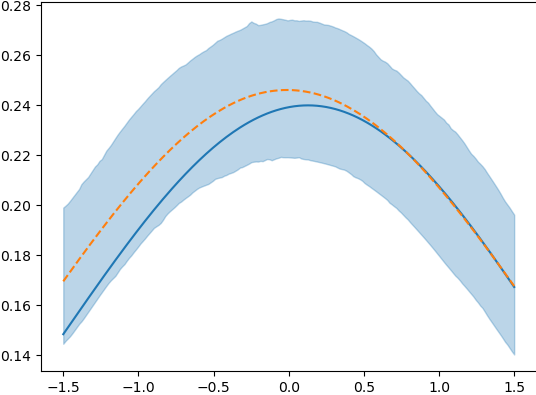}
\label{2case}
}
\caption{Performance of density estimators on $f_\epsilon,\ z$ and $g$, blue lines represent real densities(we use Monte-carlo integration to calculate the densities in (O1)), dashed line represents the point-wise average of estimators in 500 simulations, the boundaries of shaded areas respectively records the $5\%$ and $95\%$ point-wise sample quantiles in 500 simulations}
\label{Fig1}
\end{figure}

The metric we use in the first experiment is $\Vert f\Vert_{T_2} = \sqrt{\int_{[-1/h, 1/h]^d}\vert f(x)\vert^2 dx}$, which approximates $\Vert f\Vert_{L_2}$ and easily derived through Monte Carlo integrations. In table \ref{tableEst} and figure \ref{Fig1}, we observe that $\widehat{f}_\epsilon$ does not have high convergence rate. Fortunately, after convoluted with $f_\eta$, the estimators $\widehat{z}$ and $\widehat{g}$ have better performance.
\begin{table}[htbp]
  \centering
  \caption{Performance of presented estimators, we perform 500 experiments and record the sample mean (and $90\%$ quantile for crucial estimators) of each item. We do not record errors in estimating $f_X$ for (O1) and (O2) since we do not know the distribution of $X_1$ when $\epsilon_1$ is not normal}
  \begin{tabular}{l l l l  l l | l l | l l}
  \hline\hline
  Model & size & $\Vert \widehat{A} - A\Vert_2$ & $\Vert\widehat{f}_X - f_X\Vert_{T_2}$ & \multicolumn{2}{l}{$\Vert\widehat{f}_\epsilon - f_\epsilon\Vert_{T_2}$} & \multicolumn{2}{l}{$\Vert\widehat{z} - z\Vert_{T_2}$} & \multicolumn{2}{l}{$\Vert\widehat{g} - g\Vert_{T_2}$}\\
        &      & mean                           & mean                                    & mean & $90\%$ quan & mean & $90\%$ quan & mean & $90\%$ quan\\[3pt]
  \hline
  (O1)  & 500  & 0.040 & / & 0.153 & 0.230 & 0.072 & 0.112 & 0.039 & 0.059 \\
  (O1)  & 2000 & 0.021 & / & 0.124 & 0.174 & 0.053 & 0.077 & 0.030 & 0.046 \\
  (S1)  & 500  & 0.043 & 0.059 & 0.164 & 0.223 & 0.045 & 0.076 & 0.021 & 0.037\\
  (S1)  & 2000 & 0.020 & 0.046 & 0.142 & 0.171 & 0.025 & 0.042 & 0.013 & 0.021\\
  (O2)  & 5000 & 0.062 & /     & 0.173 & 0.211 & 0.080 & 0.102 & 0.050 & 0.078\\
  (S2)  & 5000 & 0.064 & 0.019 & 0.041 & 0.052 & 0.025 & 0.034 & 0.010 & 0.017\\
  \hline\hline
  \end{tabular}
  \label{tableEst}
\end{table}
\renewcommand\arraystretch{1.0}

We demonstrate the performance of algorithm \ref{alg1} in table \ref{tab2}. Two alternatives are example \ref{EXPZ} and the bootstrap algorithm used in \cite{https://doi.org/10.1111/j.1467-9892.2008.00604.x}. We do not assume knowing the state transition matrix and the distribution of state noise in algorithm \ref{alg1} but assume everything about the state space model is known in the two alternatives. The coverage probability of prediction intervals constructed by algorithm \ref{alg1} approximates the nominal coverage probability $1-\alpha$. Algorithm \ref{alg1} and example \ref{EXPZ} tend to find prediction intervals with coverage probability less than $1-\alpha$ while the bootstrap algorithm \cite{https://doi.org/10.1111/j.1467-9892.2008.00604.x} tends to find a wide prediction interval.

\begin{table}[htbp]
\centering
\caption{Performance on prediction, the result is calculated through 500 simulations, nominal coverage probability is $1-\alpha = 0.95$, F, PX and PY respectively represents filtering, predicting state vectors and predicting future observations. }
\begin{tabular}{l l l l l l l l l }
\hline\hline
Model & size & Method                    & \multicolumn{3}{l}{coverage probability} & \multicolumn{3}{l}{mean length}     \\
      &      &                           &  F       & PX      & PY           & F      & PX      & PY     \\
\hline
(O1)  & 500  & Algorithm \ref{alg1}      &  0.944   & 0.922   & 0.904        & 4.239  & 4.942   & 6.113  \\
(O1)  & 500  & Example \ref{EXPZ}        &  0.838   & 0.906   & 0.920        & 4.588  & 5.370   & 6.649  \\
(O1)  & 500  & Bootstrap                 &  /       & /       & 0.964        & /      & /       & 7.550  \\
(O1)  & 2000 & Algorithm \ref{alg1}      &  0.940   & 0.934   & 0.948        & 4.224  & 5.403   & 6.605  \\
(O1)  & 2000 & Example \ref{EXPZ}        &  0.880   & 0.904   & 0.912        & 4.588  & 5.370   & 6.649  \\
(O1)  & 2000 & Bootstrap                 &  /       & /       & 0.972        & /      & /       & 7.602  \\
(S1)  & 500  & Algorithm \ref{alg1}      &  0.960   & 0.926   & 0.910        & 3.846  & 4.776   & 5.738  \\
(S1)  & 500  & Example \ref{EXPZ}        &  0.838   & 0.896   & 0.906        & 4.588  & 5.370   & 6.649  \\
(S1)  & 500  & Bootstrap                 &  /       & /       & 0.968        & /      & /       & 7.433  \\
(S1)  & 2000 & Algorithm \ref{alg1}      &  0.960   & 0.924   & 0.930        & 3.847  & 5.050   & 6.091  \\
(S1)  & 2000 & Example \ref{EXPZ}        &  0.842   & 0.908   & 0.922        & 4.588  & 5.370   & 6.649  \\
(S1)  & 2000 & Bootstrap                 &  /       & /       & 0.980        & /      & /       & 7.519  \\
\hline\hline
\end{tabular}
\label{tab2}
\end{table}

\section{Conclusions}
Literatures in state space models focus on parametric estimation and filtering. However, if the state space model is not fully specified and maximum likelihood estimation fails, then methods like Kalman filter cannot be applied to make predictions. This paper brings a method to consistently estimate the state transition matrix and the distribution of state noises when they are unknown. In addition, this paper introduces an algorithm to construct prediction intervals of state vectors and future observations with unknown state transition matrix and distribution of state noises. In practice, fully specified state space model is too restrictive to model real life data, so methods introduced in this paper can be good alternatives to parametric estimation and filtering techniques.

\label{Conc}
\bibliographystyle{plain}
\bibliography{Draft4}

\appendix
\section{Proof of theorem \ref{thm1}}
We first prove a lemma.
\begin{lemma}
Suppose $X_i, Y_i,\ i = ...,-1,0,1,...$ are two stationary $d$ dimensional random vectors and $X_i, i = ...,-1,0,1,...$ are independent of $Y_i, i = ...,-1,0,1,...$. Define $Z_i = X_i + Y_i$ for $\forall i$, suppose the $\alpha-$mixing coefficients $\alpha_Z(k), k\geq 0\in\mathbf{Z}$ of $Z_i$ satisfy $\sum_{k = 0}^\infty \alpha_Z(k) < \infty$, then we have
\begin{equation}
\sup_{x\in\mathbf{R}^d}\mathbf{E}\vert\frac{1}{n}\sum_{i = 1}^n\exp(ix^TZ_i) - \mathcal{F}f_Z(x)\vert^2 = O(1/n)
\label{Char}
\end{equation}
\label{lemma1}
\end{lemma}

We adopt the definition in chapter 1.1 of \cite{mixing} and define the $\alpha-$mixing coefficients of $Z_i$ as
\begin{equation}
\alpha_Z(k) = \sup\{\vert Prob(U\cap V) - Prob(U)Prob(V)\vert|U\in\sigma(...,Z_0),\ V\in\sigma(Z_k,...)\}
\end{equation}
Here $\sigma(S)$ denotes the $\sigma$-field generated by $S$.

\begin{proof}[Proof of lemma \ref{lemma1}]
We define $\overline{\mathcal{H}}_k$ as the $\sigma$-field generated by $Z_k,Z_{k+1},...$ and $\underline{\mathcal{H}}_k$ as the $\sigma$-field generated by $...,Z_{k -1}, Z_k$. For any $x\in\mathbf{R}^d$, define $\overline{\mathcal{G}}_k$ as the $\sigma$-field generated by $\sin(x^TZ_i)\ \text{and}\ \cos(x^TZ_i),\ i = k,k+1,...$ and $\underline{\mathcal{G}}_k$ as the $\sigma$-field generated by $\sin(x^TZ_i)\ \text{and}\ \cos(x^TZ_i),\ i = ...,k-1,k$. Since $\sin,\cos$ are continuous function, we have $\underline{\mathcal{G}}_k\subset \underline{\mathcal{H}}_k$ and $\overline{\mathcal{G}}_k\subset \overline{\mathcal{H}}_k$, correspondingly by defining
\begin{equation}
\alpha_\mathcal{G}(k) = \sup\{\vert Prob(U\cap V) - Prob(U)Prob(V)\vert| U\in\underline{\mathcal{G}}_0,\ V\in\overline{\mathcal{G}}_k\}
\end{equation}
We have $\alpha_\mathcal{G}(k)\leq \alpha_Z(k)$ for any $k = 0,1,...$. Notice that for any $k, j$, from lemma 3 in \cite{mixing} we have
\begin{equation}
\begin{aligned}
\vert\mathbf{E}\exp(ix^TZ_k-it^TZ_j)-\mathcal{F}f_Z(x)\overline{\mathcal{F}f_Z(x)}\vert\leq \vert\mathbf{E}\cos(x^TZ_k)\cos(x^TZ_j)-\mathbf{E}\cos(x^TZ_k)\mathbf{E}\cos(x^TZ_j)\vert\\
+\vert\mathbf{E}\sin(x^TZ_k)\sin(x^TZ_j)-\mathbf{E}\sin(x^TZ_k)\mathbf{E}\sin(x^TZ_j)\vert+\vert\mathbf{E}\sin(x^TZ_k)\cos(x^TZ_j)-\mathbf{E}\sin(x^TZ_k)\mathbf{E}\cos(x^TZ_j)\vert\\
+\vert\mathbf{E}\cos(x^TZ_k)\sin(x^TZ_j)-\mathbf{E}\cos(x^TZ_k)\mathbf{E}\sin(x^TZ_j)\vert\leq 16\alpha_{\mathcal{G}}(\vert k-j\vert)\leq 16\alpha_Z(\vert k-j\vert)
\end{aligned}
\label{Cond_Mix}
\end{equation}
Therefore, we have
\begin{equation}
\begin{aligned}
\mathbf{E}\vert\frac{1}{n}\sum_{i = 1}^n\exp(ix^TZ_i) - \mathcal{F}f_Z(x)\vert^2 = \frac{1}{n^2}\sum_{i = 1}^n\sum_{j = 1}^n\mathbf{E}\left(\exp(ix^TZ_i) - \mathcal{F}f_Z(x)\right)\left(\exp(-ix^TZ_j) - \overline{\mathcal{F}f_Z(x)}\right)\\
\leq \frac{16}{n^2}\sum_{i = 1}^n\sum_{j = 1}^n\alpha_Z(\vert i - j\vert)\leq \frac{32}{n}\sum_{i = 0}^\infty\alpha_Z(i) = O(1/n)
\end{aligned}
\end{equation}
Since $x\in\mathbf{R}^d$ is arbitrarily chosen, we prove \eqref{Char}.
\end{proof}

Now we start proving theorem \ref{thm1}.

\begin{proof}[proof of theorem \ref{thm1}]
We first prove \eqref{Esti_A}. Define $Z_k = B^{-1}Y_k$ and $\tau_k = B^{-1}\eta_k$ for $k = 1,2,...,n$, we have $Z_k = X_k + \tau_k$. From assumption 1) and 3), we have $X_k = \sum_{i = 0}^\infty A^i\epsilon_{k - i}$, and correspondingly
\begin{equation}
\begin{aligned}
\mathbf{E}\Vert X_k\Vert_2^4\leq \sum_{i = 0}^\infty\sum_{j = 0}^\infty\sum_{s = 0}^\infty\sum_{l = 0}^\infty\Vert A\Vert_2^{i + j + s + l}\mathbf{E}\Vert\epsilon_{k - i}\Vert_2\Vert\epsilon_{k - j}\Vert_2\Vert\epsilon_{k - s}\Vert_2\Vert\epsilon_{k - l}\Vert_2\leq\mathbf{E}\Vert\epsilon_1\Vert_2^4\times \left(\sum_{k = 0}^\infty \Vert A\Vert_2^k\right)^4<\infty
\end{aligned}
\end{equation}
This implies $\Sigma_X = \mathbf{E}X_1X_1^T$ exists and $\Sigma_X = \sum_{i = 0}^\infty A^i\Sigma A^{iT}$ is positive definite. Correspondingly
\begin{equation}
\begin{aligned}
Z_kZ_{k-2}^T = A^2X_{k-2}X_{k-2}^T + A\epsilon_{k-1}X_{k-2}^T + \epsilon_kX_{k-2}^T + \tau_kX^T_{k-2} + X_k\tau_{k-2}^T + \tau_k\tau^T_{k - 2}\Rightarrow \mathbf{E}Z_kZ_{k-2}^T = A^2\Sigma_X\\
Z_kZ_{k - 1}^T = AX_{k - 1}X_{k - 1}^T + \epsilon_kX^T_{k - 1} + \tau_kX^T_{k - 1} + X_k\tau^T_{k - 1} + \tau_k\tau_{k - 1}^T\Rightarrow \mathbf{E}Z_kZ_{k - 1}^T = A\Sigma_X
\label{Y_k}
\end{aligned}
\end{equation}
For a matrix $A$, we define its Frobenius norm as $\Vert A\Vert_F = \sqrt{Tr(AA^T)}$, since
\begin{equation}
\begin{aligned}
\mathbf{E}\Vert\frac{1}{n - 2}\sum_{k = 1}^{n - 2}X_kX_k^T - \Sigma_X\Vert_F^2 = \frac{1}{(n - 2)^2}\sum_{k = 1}^{n - 2}\sum_{j = 1}^{n - 2}Tr\left(\mathbf{E}(X_kX_k^T - \Sigma_X)(X_jX_j^T - \Sigma_X)\right)\\
=\frac{1}{(n - 2)^2}\sum_{k = 1}^{n - 2}\sum_{j = 1}^{n - 2}\left(\mathbf{E}(X_k^TX_j)^2 - Tr(\Sigma_X^2)\right)\\
=\frac{\mathbf{E}\Vert X_1\Vert_2^4 - Tr(\Sigma_X^2)}{n - 2} + \frac{2}{(n - 2)^2}\sum_{1\leq k<j\leq n-2}\left(\mathbf{E}\left(\sum_{s=0}^{j-k-1}X_k^TA^s\epsilon_{j-s}+X_k^TA^{j-k}X_k\right)^2- Tr\Sigma_X^2\right)\\
= \frac{\mathbf{E}\Vert X_1\Vert_2^4 - Tr(\Sigma_X^2)}{n - 2} +  \frac{2}{(n - 2)^2}\sum_{1\leq k < j \leq n - 2}\left(Tr(\sum_{s = 0}^{j - k - 1}A^s\Sigma A^{sT}\Sigma_X - \Sigma^2_X) + \mathbf{E}(X_k^TA^{j - k}X_k)^2\right)\\
\leq \frac{\mathbf{E}\Vert X_1\Vert_2^4 - Tr(\Sigma_X^2)}{n - 2} + \frac{2}{(n - 2)^2}\sum_{1\leq k < j\leq n- 2}\left(\Vert\sum_{s= 0}^{j - k - 1}A^s\Sigma A^{sT} - \Sigma_X\Vert_F\times\Vert\Sigma_X\Vert_F + \Vert A\Vert_2^{2j-2k}\mathbf{E}\Vert X_k\Vert_2^4\right)\\
\leq \frac{\mathbf{E}\Vert X_1\Vert_2^4 - Tr(\Sigma_X^2)}{n - 2} +\frac{2\sqrt{d}\Vert\Sigma_X\Vert_F\Vert\Sigma\Vert_2}{n - 2}\sum_{k = 1}^{n -3}\frac{\Vert A\Vert_2^{2k}}{1 - \Vert A\Vert_2^2} +\frac{2\mathbf{E}\Vert X_k\Vert_2^4}{n - 2}\sum_{k = 1}^{n - 3}\Vert A\Vert_2^{2k} = O(1/n)
\end{aligned}
\end{equation}
Similarly we have
\begin{equation}
\begin{aligned}
\mathbf{E}\Vert\frac{1}{n-2}\sum_{k=1}^{n-2}\epsilon_{k+1}X_k^T\Vert_F^2= \frac{1}{(n-2)^2}\sum_{k=1}^{n-2}\sum_{j=1}^{n-2}\mathbf{E}X_k^TX_j\epsilon_{j+1}^T\epsilon_{k+1}=\frac{1}{n-2}\mathbf{E}\Vert X_1\Vert_2^2\mathbf{E}\Vert\epsilon_1\Vert_2^2 = O(1/n)\\
\mathbf{E}\Vert\frac{1}{n-2}\sum_{k = 1}^{n - 2}\epsilon_{k+2}X_k^T\Vert^2_F = \frac{1}{(n - 2)^2}\sum_{k = 1}^{n-2}\sum_{j = 1}^{n - 2}\mathbf{E}X_k^TX_{j}\epsilon_{k+2}^T\epsilon_{j +2 } = \frac{1}{n-2}\mathbf{E}\Vert X_1\Vert_2^2\mathbf{E}\Vert\epsilon_1\Vert_2^2 = O(1/n)\\
\mathbf{E}\Vert \frac{1}{n-2}\sum_{k=1}^{n-2}X_k\tau_{k+2}^T\Vert_F^2= \frac{1}{(n-2)^2}\sum_{k=1}^{n-1}\sum_{j=1}^{n-2}\mathbf{E}X_k^TX_j\tau_{j+2}^T\tau_{k+2}=\frac{1}{n-2}\mathbf{E}\Vert X_1\Vert_2^2\mathbf{E}\Vert\tau_1\Vert_2^2 = O(1/n)\\
\mathbf{E}\Vert\frac{1}{n - 2}\sum_{k = 1}^{n - 2}X_{k + 2}\tau_k\Vert_F = \frac{1}{(n - 2)^2}\sum_{k = 1}^{n - 2}\sum_{j = 1}^{n - 2}\mathbf{E}X_{k+2}^TX_{j+2}\tau_k^T\tau_j = \frac{1}{n - 2}\mathbf{E}\Vert X_1\Vert_2^2\mathbf{E}\Vert\tau_1\Vert_2^2 = O(1/n)\\
\mathbf{E}\Vert\frac{1}{n-2}\sum_{k=1}^{n-2}\epsilon_{k+2}\tau_k^T\Vert_F^2 = \frac{1}{(n-2)^2}\sum_{k=1}^{n-2}\sum_{j=1}^{n-2}\mathbf{E}\epsilon_{j+2}^T\epsilon_{k+2}\tau_{k}^T\tau_j=\frac{1}{n-2}\mathbf{E}\Vert\epsilon_1\Vert_2^2\mathbf{E}\Vert\tau_1\Vert_2^2 = O(1 / n)\\
\mathbf{E}\Vert\frac{1}{n-2}\sum_{k=1}^{n-2}\tau_{k+2}\tau_k^T\Vert_F^2 = \frac{1}{(n-2)^2}\sum_{k=1}^{n-2}\sum_{j=1}^{n-2}\mathbf{E}\tau_k^T\tau_j\tau_{j+2}^T\tau_{k+2}=\frac{1}{n-2}(\mathbf{E}\Vert\tau_1\Vert_2^2)^2 = O(1/n)
\end{aligned}
\label{Df}
\end{equation}
Therefore, from \eqref{Y_k} to \eqref{Df}, we know that $\Vert\frac{1}{n - 2}\sum_{k = 3}^n Z_kZ_{k-2}^T - A^2\Sigma_X\Vert_2 = O_p(1/\sqrt{n})$ and similarly $\Vert\frac{1}{n - 2}\sum_{k = 3}^n Z_{k - 1}Z_{k - 2}^T - A\Sigma_X\Vert_2 = O_p(1/\sqrt{n})$. Apply corollary 6.3.4 and (5.8.4) in \cite{Horn:2012:MA:2422911}, since $A\Sigma_X$ is non-singular, for sufficiently large $n$, with large probability $\frac{1}{n - 2}\sum_{k = 3}^n Z_{k - 1}Z_{k - 2}^T $ is invertible, and we have $\Vert(\frac{1}{n - 2}\sum_{k = 3}^n Z_{k - 1}Z_{k - 2}^T)^{-1} - (A\Sigma_X)^{-1}\Vert_2 = O_p(1/\sqrt{n})$ as well. Therefore, we have
\begin{equation}
\begin{aligned}
\Vert\widehat{A} - A\Vert_2 = \Vert\left(\frac{1}{n - 2}\sum_{k = 3}^n Z_kZ_{k-2}^T\right)\left(\frac{1}{n - 2}\sum_{k = 3}^n Z_{k - 1}Z_{k - 2}^T \right)^{+} - (A^2\Sigma_X)(A\Sigma_X)^{-1}\Vert_2\\
\leq \Vert\frac{1}{n - 2}\sum_{k = 3}^n Z_kZ_{k-2}^T\Vert_2\times\Vert \left(\frac{1}{n - 2}\sum_{k = 3}^n Z_{k - 1}Z_{k - 2}^T \right)^{+} - (A\Sigma_X)^{-1}\Vert_2 + \Vert \left(\frac{1}{n - 2}\sum_{k = 3}^n Z_kZ_{k-2}^T\right) - A^2\Sigma_X\Vert_2\times\Vert(A\Sigma_X)^{-1}\Vert_2
\end{aligned}
\end{equation}
Which is of $O_p(1/\sqrt{n})$ and we prove the first result.

From Plancherel theorem and assumption 3), we have
\begin{equation}
\begin{aligned}
\mathbf{E}\widetilde{f}_X - f_X= \mathcal{F}^{-1}\left((\mathcal{F}K_h - 1)\mathcal{F}f_X\right)\\
\Rightarrow \Vert \mathbf{E}\widetilde{f}_X - f_X\Vert_{L_2} = \frac{\Vert (\mathcal{F}K_h - 1)\mathcal{F}f_X\Vert_{L_2}}{(2\pi)^{d/2}}
\leq \frac{1}{(2\pi)^{d/2}}\sqrt{\sum_{i = 1}^d \int_{\vert x_i\vert > \frac{1}{h}}\vert\mathcal{F}f_X(x)\vert^2 dx}
\end{aligned}
\label{Bias1}
\end{equation}
If assumption 4.a) happens, we have for any $i = 1,2,...,d$,
\begin{equation}
\begin{aligned}
\int_{\vert x_i\vert > \frac{1}{h}}\vert\mathcal{F}f_X(x)\vert^2 dx\leq \frac{h^b}{(1 + h^2)^{b / 2}}\int_{\mathbf{R}^d}\vert\mathcal{F}f_X(x)\vert^2(1 + x_i^2)^{b / 2}dx
\Rightarrow\Vert\mathbf{E}\widetilde{f}_X - f_X\Vert_{L_2} = O(h^{b / 2})
\end{aligned}
\label{Biasa}
\end{equation}
If assumption 4.b) happens, we have for any $i = 1,2,...,d$,
\begin{equation}
\begin{aligned}
\int_{\vert x_i\vert > \frac{1}{h}}\vert\mathcal{F}f_X(x)\vert^2dx\leq \exp\left(-\frac{r}{h^b}\right)\int_{\mathbf{R}^d}\vert\mathcal{F}f_X(x)\vert^2\exp(r\vert x_i\vert^b)dx\Rightarrow \Vert \mathbf{E}\widetilde{f}_X - f_X\Vert_{L_2} = O\left(\exp\left(-\frac{r}{2h^b}\right)\right)
\end{aligned}
\label{Bias2}
\end{equation}
Notice that
\begin{equation}
\left[
\begin{matrix}
Z_{n+1}\\
\tau_{n+1}
\end{matrix}
\right]=\left[
\begin{matrix}
A & -A\\
0 & 0
\end{matrix}
\right]
\left[
\begin{matrix}
Z_{n}\\
\tau_{n}
\end{matrix}
\right]+\left[
\begin{matrix}
I & I\\
0 & I
\end{matrix}
\right]\left[
\begin{matrix}
\epsilon_{n+1}\\
\tau_{n+1}
\end{matrix}
\right]
\end{equation}
is an ARMA model and $A$ is non-singular with $\Vert A\Vert_2<1$, so the polynomial
\begin{equation}
f(x)=\det\left(
\left[
\begin{matrix}
I & I\\
0 & I
\end{matrix}
\right]-\left[
\begin{matrix}
A & -A\\
0 & 0
\end{matrix}
\right]x
\right)=det(I-Ax)
\end{equation}
has all roots with absolute value greater than 1. According to theorem 1 in \cite{MOKKADEM1988309}, $\left[Z_i,\tau_i\right]^T$ is a mixing process and $\exists 0<s<1$ such that the $\alpha-$mixing coefficients $\alpha_{[Z_i,\tau_i]^T}(k)=O(s^{\vert k\vert})$ for any $k\in \mathbf{Z}$. Combine with lemma \ref{lemma1} and assumption 3), we have
\begin{equation}
\begin{aligned}
\widetilde{f}_X(x) - \mathbf{E}\widetilde{f}_X(x) = \frac{1}{(2\pi)^d}\int_{\mathbf{R}^d}\frac{\mathcal{F}K_h(y)}{\mathcal{F}f_\eta(B^{-1T}y)}\left(\frac{1}{n}\sum_{j = 1}^n\exp(iy^TZ_j) - \mathcal{F}f_Z(y)\right)\exp(-iy^Tx)dy\\
\Rightarrow \mathbf{E}\Vert\widetilde{f}_X(x) - \mathbf{E}\widetilde{f}_X(x)\Vert^2_{L_2} = \frac{1}{(2\pi)^d}\int_{\mathbf{R}^d}\vert\frac{\mathcal{F}K_h(y)}{\mathcal{F}f_\eta(B^{-1T}y)}\vert^2\mathbf{E}\vert \frac{1}{n}\sum_{j = 1}^n\exp(iy^TZ_j) - \mathcal{F}f_Z(y)\vert^2dy\\
= O\left(\frac{1}{n}\int_{\mathbf{R}^d}\vert\frac{\mathcal{F}K_h(y)}{\mathcal{F}f_\eta(B^{-1T}y)}\vert^2dy\right)
= O\left(\frac{1}{n}\int_{\Vert y\Vert_\infty\leq \frac{a}{h}}\frac{1}{\vert\mathcal{F}f_\eta(B^{-1T}y)\vert^2}dy\right)
\end{aligned}
\label{Var1}
\end{equation}
Since $B$ is non-singular, by change of variable theorem,
\begin{equation}
\begin{aligned}
\int_{\Vert y\Vert_\infty\leq \frac{a}{h}}\frac{1}{\vert\mathcal{F}f_\eta(B^{-1T}y)\vert^2}dy\leq \int_{\Vert B^{-1T}y\Vert_\infty\leq \frac{a} {hc_{B^{-1T}}}}\frac{1}{\vert\mathcal{F}f_\eta(B^{-1T}y)\vert^2}dy = \int_{\Vert y\Vert_\infty\leq \frac{a}{hc_{B^{-1T}}}}\frac{\vert\det B\vert}{\vert \mathcal{F}f_\eta(y)\vert^2}dy
\end{aligned}
\label{Equi1}
\end{equation}
If assumption 4.a) happens, according to lemma 1. in \cite{AIHPB_2013__49_2_569_0} and \eqref{Biasa},
\begin{equation}
\begin{aligned}
\int_{\Vert y\Vert_\infty\leq \frac{a}{hc_{B^{-1T}}}}\frac{1}{\vert \mathcal{F}f_\eta(y)\vert^2}dy = O\left( \int_{\Vert y\Vert_\infty\leq \frac{a}{hc_{B^{-1T}}}}\prod_{j = 1}^d (y_j^2 + 1)^{\beta}dy\right) = O\left(\frac{1}{h^{2d\beta + d}}\right)\\
\Rightarrow \mathbf{E}\Vert\widehat{f}_X -f_X\Vert_{L_2}^2\leq \mathbf{E}\Vert\widetilde{f}_X -f_X\Vert^2_{L_2}
\leq 2\mathbf{E}\Vert\widetilde{f}_X - \mathbf{E}\widetilde{f}_X\Vert^2_{L_2} + 2\Vert \mathbf{E}\widetilde{f}_X - f_X\Vert^2_{L_2}
= O\left(\frac{1}{nh^{2d\beta + d}} + h^b\right)
\end{aligned}
\end{equation}
which tends to $0$ as $n\to\infty$. If assumption 4.b) happens, similarly we have
\begin{equation}
\begin{aligned}
\int_{\Vert y\Vert_\infty\leq \frac{a}{hc_{B^{-1T}}}}\frac{1}{\vert \mathcal{F}f_\eta(y)\vert^2}dy = O\left(h^{d\beta - d}\exp\left(2d\gamma\frac{a^\beta}{h^\beta c^\beta_{B^{-1T}}}\right)\right)\\
\Rightarrow \mathbf{E}\Vert\widehat{f}_X -f_X\Vert^2_{L_2} = O\left(\exp\left(2d\gamma\frac{a^\beta}{h^\beta c^\beta_{B^{-1T}}}-\log(n) + (d\beta - d)\log(h)\right) + \exp\left(-\frac{r}{h^b}\right)\right)
\end{aligned}
\label{Var2}
\end{equation}
$h = o(1)\Rightarrow \exp\left(-\frac{r}{h^b}\right) = o(1)$. On the other hand, let $z = h\log^{1/\beta}(n)\to\infty$, for sufficiently large $n$,
\begin{equation}
2d\gamma\frac{a^\beta}{h^\beta c^\beta_{B^{-1T}}}-\log(n) + (d\beta - d)\log(h) = (\frac{2d\gamma a^\beta}{ z^\beta c^\beta_{B^{-1T}}} - 1)\log(n) + (d - d\beta)\log(\frac{1}{h})\leq -\frac{1}{2}\log(n) + \frac{d}{\beta}\log(\log(n))
\label{Zs2}
\end{equation}
which tends to $-\infty$ as $n\to\infty$ and we prove 2.

We define
\begin{equation}
\widetilde{f}_\epsilon(x)^\dagger = \frac{1}{(2\pi)^d\times (n - 1)}\sum_{j = 1}^{n - 1}\int_{\mathbf{R}^d}\frac{\mathcal{F}K_h(y)}{\mathcal{F}f_\eta(B^{-1T}y)\mathcal{F}f_\eta(-B^{-1T}A^Ty)}\exp(iy^TZ_{j+1} - it^TAZ_j - iy^Tx)dy
\end{equation}
Similar with \eqref{Bias1} to \eqref{Bias2} and \eqref{Var1} to \eqref{Var2}, if assumption 4.a) happens,
\begin{equation}
\begin{aligned}
\mathbf{E}\widetilde{f}_\epsilon(x)^\dagger = \mathcal{F}\left(\mathcal{F}K_h\times\mathcal{F}f_\epsilon\right)(x)\Rightarrow \Vert\mathbf{E}\widetilde{f}_\epsilon(x)^\dagger - f_\epsilon(x)\Vert_{L_2} = O\left(h^{b/2}\right)
\end{aligned}
\end{equation}
Since $Z_j - AZ_{j - 1}$ is a continuous function of $Z_j, Z_{j - 1}$, it is a measurable function in the $\sigma$-field generated by $...,Z_{j-1},Z_j$ or $Z_{j-1}, Z_{j},...$, this implies for any $k > 1$, the $\alpha$-mixing coefficients satisfy $\alpha_{Z_j - AZ_{j - 1}}(k)\leq \alpha_Z(k - 1)\Rightarrow \sum_{k = 0}^\infty \alpha_{Z_j - AZ_{j - 1}}(k) < \infty$, correspondingly from lemma \ref{lemma1} and Cauchy inequality

\begin{equation}
\begin{aligned}
(2\pi)^d\times \mathbf{E}\Vert\widetilde{f}_\epsilon(x)^\dagger - \mathbf{E}\widetilde{f}_\epsilon(x)^\dagger\Vert^2_{L_2}\\
=\int_{\mathbf{R}^d}\frac{\vert \mathcal{F}K_h(y)\vert^2\times \mathbf{E}\vert\sum_{j = 1}^{n - 1}\left(\exp(iy^T(Z_{j + 1} -AZ_j)) -\mathcal{F}f_\eta(B^{-1T}y)\mathcal{F}f_\epsilon(y)\mathcal{F}f_\eta(-B^{-1T}A^Ty)\right)\vert^2}{(n - 1)^2\times\vert \mathcal{F}f_\eta(B^{-1T}y)\vert^2\times \vert \mathcal{F}f_\eta(-B^{-1T}A^Ty)\vert^2}dy\\
= O\left(\frac{1}{n}\int_{\mathbf{R}^d}\frac{\vert \mathcal{F}K_h(y)\vert^2}{\vert \mathcal{F}f_\eta(B^{-1T}y)\vert^2 \times \vert \mathcal{F}f_\eta(-B^{-1T}A^Ty)\vert^2}dy\right)\\
= O\left(\frac{1}{n}\sqrt{\int_{\mathbf{R}^d}\frac{\vert \mathcal{F}K_h(y)\vert^2}{\vert \mathcal{F}f_\eta(B^{-1T}y)\vert^4}dy}\times\sqrt{\int_{\mathbf{R}^d}\frac{\vert \mathcal{F}K_h(y)\vert^2}{\vert \mathcal{F}f_\eta(-B^{-1T}A^Ty)\vert^4}dy}\right)
= O\left(\frac{1}{nh^{4\beta d + d}}\right)\\
\Rightarrow \mathbf{E}\Vert \widetilde{f}_\epsilon(x)^\dagger - f_\epsilon(x)\Vert^2_{L_2} = O(h^b + \frac{1}{nh^{4d\beta + d}})
\end{aligned}
\label{FKH}
\end{equation}

If assumption 4.b) happens, we define $c_1 = c_{B^{-1T}}$ and $c_2 = c_{B^{-1T}A^T}$,
\begin{equation}
\begin{aligned}
\mathbf{E}\Vert\widetilde{f}_\epsilon(x)^\dagger - f_\epsilon(x)\Vert^2_{L_2} = O\left(\exp\left(-\frac{r}{h^b}\right) + \frac{1}{n}\sqrt{\int_{\mathbf{R}^d}\frac{\vert \mathcal{F}K_h(y)\vert^2}{\vert \mathcal{F}f_\eta(B^{-1T}y)\vert^4}dy}\times\sqrt{\int_{\mathbf{R}^d}\frac{\vert \mathcal{F}K_h(y)\vert^2}{\vert \mathcal{F}f_\eta(-B^{-1T}A^Ty)\vert^4}dy}\right)\\
= O\left(\exp\left(-\frac{r}{h^b}\right) + \exp(2d\gamma\frac{a^\beta}{h^\beta c_1^\beta} + 2d\gamma\frac{a^\beta}{h^\beta c_2^\beta} + (d\beta - d)\log(h)-\log(n))\right)
\end{aligned}
\label{Eq33}
\end{equation}
Similar with \eqref{Zs2}, we have $\mathbf{E}\Vert\widetilde{f}_\epsilon(x)^\dagger - f_\epsilon(x)\Vert^2_{L_2} \to 0$ as $n\to\infty$.

Since for any $x,y\in\mathbf{R}^d$,
\begin{equation}
\vert\mathcal{F}f_\eta(x) - \mathcal{F}f_\eta(y)\vert\leq \int_{\mathbf{R}^d}\vert\exp(it^Tx) - \exp(it^Ty)\vert f_\eta(t)dt\leq \int_{\mathbf{R}^d}\vert(x -y)^Tt\vert f_\eta(t)dt\leq \Vert x - y \Vert_2\int_{\mathbf{R}^d}\Vert t\Vert_2 f_\eta(t)dt
\end{equation}
from Cauchy's inequality and theorem 8.22 in \cite{Folland:1706460}, by defining $C = (\int_{\mathbf{R}^d}\Vert t\Vert_2 f_\eta(t)dt)^2<\infty$,
\begin{equation}
\begin{aligned}
(2\pi)^d\Vert\widetilde{f}_\epsilon(x) - \widetilde{f}_\epsilon^\dagger(x)\Vert^2_{L_2}\\
\leq2\int_{\mathbf{R}^d}\frac{\vert\mathcal{F}K_h(y)\vert^2\times\vert\sum_{j = 1}^{n - 1}\left(\exp(iy^TZ_{j+1} - iy^T\widehat{A}Z_j) -\exp(iy^TZ_{j+1} - iy^TAZ_j)\right)\vert^2}{(n - 1)^2\times\vert\mathcal{F}f_\eta(B^{-1T}y)\vert^2\times\vert\mathcal{F}f_\eta(-B^{-1T}\widehat{A}^Ty)\vert^2}dy\\
+2\int_{\mathbf{R}^d}\frac{\vert\mathcal{F}K_h(y)\vert^2\times \vert\sum_{j = 1}^{n - 1}\exp(iy^TZ_{j+1} - iy^TAZ_j)\vert^2}{(n - 1)^2\times\vert\mathcal{F}f_\eta(B^{-1T}y)\vert^2}\times \vert\frac{1}{\mathcal{F}f_\eta(-B^{-1T}\widehat{A}^Ty)} - \frac{1}{\mathcal{F}f_\eta(-B^{-1T}A^Ty)}\vert^2dy\\
\leq \frac{8}{n - 1}\int_{\mathbf{R}^d}\frac{\vert\mathcal{F}K_h(y)\vert^2\sum_{j = 1}^{n - 1}\sin^2\left(\frac{y^T(A - \widehat{A})Z_j}{2}\right)}{\vert\mathcal{F}f_\eta(B^{-1T}y)\vert^2\times\vert\mathcal{F}f_\eta(-B^{-1T}\widehat{A}^Ty)\vert^2}dy\\
+ 2\int_{\mathbf{R}^d}\frac{\vert\mathcal{F}K_h(y)\vert^2\times \vert \sum_{j = 1}^{n - 1}\exp(iy^TZ_{j+1} - iy^TAZ_j)\vert^2}{(n - 1)^2\times\vert\mathcal{F}f_\eta(B^{-1T}y)\vert^2\times \vert\mathcal{F}f_\eta(-B^{-1T}A^Ty)\vert^2}\times\frac{\vert\mathcal{F}f_\eta(-B^{-1T}\widehat{A}^Ty) -\mathcal{F}f_\eta(-B^{-1T}A^Ty)\vert^2}{\vert\mathcal{F}f_\eta(-B^{-1T}\widehat{A}^Ty)\vert^2}dy\\
\leq \left(\int_{\mathbf{R}^d}\frac{\vert\mathcal{F}K_h(y)\vert^2}{\vert\mathcal{F}f_\eta(B^{-1T}y)\vert^2\times \vert\mathcal{F}f_\eta(-B^{-1T}A^Ty)\vert^2}dy\right)\times\frac{2}{(n-1)}\max_{\Vert y\Vert_\infty\leq\frac{a}{h}}\frac{\vert\mathcal{F}f_\eta(-B^{-1T}A^Ty)\vert^2}{\vert\mathcal{F}f_\eta(-B^{-1T}\widehat{A}^Ty)\vert^2}\times\sum_{j = 1}^{n - 1}\vert y^T(A - \widehat{A})Z_j\vert^2\\
+ 2\int_{\mathbf{R}^d}\frac{\vert\mathcal{F}K_h(y)\vert^2\times \vert\sum_{j = 1}^{n - 1}\exp(iy^TZ_{j+1} - iy^TAZ_j)\vert^2}{(n - 1)^2\times\vert\mathcal{F}f_\eta(B^{-1T}y)\vert^2\times \vert\mathcal{F}f_\eta(-B^{-1T}A^Ty)\vert^2}dy\times\max_{\Vert y\Vert_\infty\leq \frac{a}{h}}\frac{C\Vert B^{-1T}(\widehat{A} - A)^Ty\Vert^2_2}{\vert\mathcal{F}f_\eta(-B^{-1T}\widehat{A}^Ty)\vert^2}
\end{aligned}
\label{Transation}
\end{equation}
Since $\mathbf{E}\frac{1}{n - 1}\sum_{j = 1}^{n - 1}\Vert Z_j\Vert_2^2 = \mathbf{E}\Vert Z_1\Vert_2^2 < \infty\Rightarrow \frac{1}{n - 1}\sum_{j = 1}^{n - 1}\Vert Z_j\Vert_2^2  = O_p(1)$, if 4.a) is satisfied
\begin{equation}
\begin{aligned}
\max_{\Vert y\Vert_\infty\leq \frac{a}{h}}\frac{\vert\mathcal{F}f_\eta(-B^{-1T}A^Ty)\vert^2}{\vert\mathcal{F}f_\eta(-B^{-1T}\widehat{A}^Ty)\vert^2}\leq 2 + \max_{\Vert y\Vert_\infty\leq \frac{a}{h}}\frac{2C\Vert B^{-1}\Vert_2^2\Vert A - \widehat{A}\Vert_2^2\Vert y\Vert_2^2}{\vert\mathcal{F}f_\eta(-B^{-1T}\widehat{A}^Ty)\vert^2}\\
\leq 2 + \frac{2C\Vert B^{-1}\Vert_2^2\Vert A - \widehat{A}\Vert_2^2 da^2}{c^2h^2} (\Vert B^{-1}\Vert_2^2\Vert\widehat{A}\Vert_2^2\frac{da^2}{h^2} + 1)^{d\beta} = O_p(1 + \frac{1}{nh^{2d\beta + 2}})
\end{aligned}
\end{equation}
combine with \eqref{FKH} and \eqref{Transation}, we have
\begin{equation}
\begin{aligned}
\Vert\widetilde{f}_\epsilon(x) - \widetilde{f}_\epsilon^\dagger(x)\Vert^2_{L_2} = O_p\left(\frac{1}{nh^{4d\beta+d+2}} + \frac{1}{n^2h^{6d\beta + 4 + d}} + \frac{1}{nh^{2d\beta + 2}}\right) = o_p(1)
\end{aligned}
\end{equation}

On the other hand, if 4.b) is satisfied, by letting $s = h\log^{1/\beta}(n)\to\infty$, if $\Vert\widehat{A}\Vert_2\leq 2^{1/\beta}\Vert A\Vert_2$, for sufficient large $n$, $\frac{4d\gamma\Vert B^{-1}\Vert_2^\beta\Vert A\Vert_2^\beta\times d^{\beta / 2}a^\beta}{s^\beta} < 1/4$, correspondingly
\begin{equation}
\begin{aligned}
\max_{\Vert t\Vert_\infty\leq a/h}\frac{\vert\mathcal{F}f_\eta(-B^{-1T}A^Tt)\vert^2}{\vert\mathcal{F}f_\eta(-B^{-1T}\widehat{A}^Tt)\vert^2}
\leq 2 + \frac{2C\Vert B^{-1}\Vert_2^2\Vert A - \widehat{A}\Vert_2^2 da^2}{c^2h^2}\exp\left(2d\gamma\Vert B^{-1}\Vert^\beta_2\Vert\widehat{A}\Vert^\beta_2\times \log(n)\frac{d^{\beta/2}a^\beta}{s^\beta}\right)\\
\leq 2 + C^{'}n^{1/4}\frac{\Vert\widehat{A} - A\Vert_2^2}{h^2}
\end{aligned}
\end{equation}
combine with \eqref{Eq33}, since $\Vert\widehat{A} - A\Vert_2\to_p 0$ and \eqref{Zs2}, we have
\begin{equation}
\begin{aligned}
\Vert\widetilde{f}_\epsilon(x) - \widetilde{f}_\epsilon^\dagger(x)\Vert^2_{L_2}\\
= O_p\left(\exp(2d\gamma\left(\frac{a}{hc_1}\right)^\beta + 2d\gamma\left(\frac{a}{hc_2}\right)^\beta + (d\beta - d)\log(h))\times (1 + \frac{1}{n^{3/4}h^2})\times \frac{1}{nh^2}+\frac{1}{n^{3/4}h^2}\right)\\
= O_p\left(\frac{\log^{(2 + d)/\beta}(n)}{\sqrt{n}}\right)
\end{aligned}
\end{equation}
Thus, we have
\begin{equation}
\Vert \widehat{f}_\epsilon(x) - f_\epsilon(x)\Vert^2_{L_2}\leq 2\Vert\widetilde{f}(x) - \widetilde{f}^\dagger(x)\Vert^2_{L_2}+ 2\Vert\widetilde{f}^\dagger(x) - f(x)\Vert^2_{L_2} = o_p(1)
\end{equation}
\end{proof}
\section{Proofs for theorems in chapter \ref{cp5}}
\begin{proof}[Proof of lemma \ref{lemma2}]
For the first result, notice that
\begin{equation}
\begin{aligned}
Prob\left(\Vert X_n - B^{-1}Y_n\Vert_\infty\leq x\right) = Prob\left(\Vert B^{-1}\eta_n\Vert_\infty\leq x\right) = \int_{\Vert B^{-1}y\Vert_\infty\leq x}f_\eta(y)dy
\end{aligned}
\end{equation}
For the second result, we have
\begin{equation}
\begin{aligned}
Prob\left(\Vert X_{n + 1} - AB^{-1}Y_n\Vert_\infty\leq x\right) = \mathbf{E}Prob\left(\Vert \epsilon_{n+1} - AB^{-1}\eta_n\Vert_\infty\leq x|\eta_n\right)\\
= \mathbf{E}\int_{\Vert z\Vert_\infty\leq x}f_\epsilon(z + AB^{-1}\eta_n)dz = \int_{\mathbf{R}^d}\left(\int_{\Vert z\Vert_\infty\leq x}f_\epsilon(z + AB^{-1}y)dz\right)f_\eta(y)dy
\end{aligned}
\end{equation}
For the third result, from the change of variable theorem
\begin{equation}
\begin{aligned}
Prob\left(\Vert Y_{n + 1} - BAB^{-1}Y_n\Vert_\infty\leq x\right) = \mathbf{E}Prob\left(\Vert\eta_{n+1} + B\epsilon_{n+1}-BAB^{-1}\eta_n\Vert_\infty\leq x|\eta_n, \epsilon_{n+1}\right)\\
=\mathbf{E}\int_{\Vert z\Vert_\infty\leq x}f_\eta(z + BAB^{-1}\eta_n - B\epsilon_{n+1})dz =\int_{\mathbf{R}^d\times\mathbf{R}^d}\left(\int_{\Vert z\Vert_\infty\leq x}f_\eta(z + BAB^{-1}w - By)dz\right)f_\eta(w)f_\epsilon(y)dwdy\\
=\frac{1}{\vert\det(B)\vert}\int_{\mathbf{R}^d\times\mathbf{R}^d}\left(\int_{\Vert y\Vert_\infty\leq x}f_\epsilon(B^{-1}y + AB^{-1}w - B^{-1}z)dy\right)f_\eta(w)f_\eta(z)dwdz
\end{aligned}
\end{equation}
\end{proof}

\begin{proof}[Proof of theorem \ref{MainThm}]
Notice that for any given $x_0 > 0$
\begin{equation}
\begin{aligned}
\sup_{0<x\leq x_0}\vert N(x) - Prob\left(\Vert X_{n+1} - AB^{-1}Y_n\Vert_\infty\leq x\right)\vert\\
\leq\sup_{0<x\leq x_0, y\in\mathbf{R}^d}\int_{\Vert z\Vert_\infty\leq x}\vert\widehat{f}_\epsilon(z + \widehat{A}B^{-1}y)dz -f_\epsilon(z + \widehat{A}B^{-1}y)\vert dz\\
+ \sup_{0<x\leq x_0, y\in\mathbf{R}^d}\int_{\Vert z\Vert_\infty\leq x}\vert f_\epsilon(z + \widehat{A}B^{-1}y) - f_\epsilon(z + AB^{-1}y)\vert dz\\
\leq (2x_0)^{d/2}\Vert\widehat{f}_\epsilon - f_\epsilon\Vert_{L_2} + (2x_0)^d\sup_{\Vert z\Vert_\infty\leq x_0, y\in\mathbf{R}^d}\vert f_\epsilon(z + \widehat{A}B^{-1}y) - f_\epsilon(z + AB^{-1}y)\vert
\end{aligned}
\label{Decop1}
\end{equation}
From remark \ref{remark1}, for any $\delta > 0$, there exists $M > 0$ such that $f_\epsilon(x)<\delta$ for any $\Vert x\Vert_2 > M$, from theorem \ref{thm1}, there exists $C_\delta > 0$ such that $\Vert\widehat{A} - A\Vert_2\leq C_\delta/\sqrt{n}$ with probability at least $1-\delta$, if this happens, for sufficiently large $n$
\begin{equation}
\begin{aligned}
\sup_{\Vert z\Vert_\infty\leq x_0, \Vert y\Vert_2\leq n^{1/4}}\vert f_\epsilon(z + AB^{-1}y) - f_\epsilon(z + \widehat{A}B^{-1}y)\vert\leq \sup_{\Vert x - y\Vert_2\leq C_\delta\Vert B^{-1}\Vert_2/n^{1/4}}\vert f_\epsilon(x) - f_\epsilon(y)\vert < \delta
\end{aligned}
\end{equation}
and for $\forall z,y$ such that $\Vert z\Vert_\infty\leq x_0$ and $\Vert y\Vert_2>n^{1/4}$, we have $\Vert z + AB^{-1}y\Vert_2\geq \frac{\Vert y\Vert_\infty}{C_{AB^{-1}}}-\Vert z\Vert_2\geq \frac{n^{1/4}}{C_{AB^{-1}}\sqrt{d}} - \sqrt{d}x_0>M$ and $\Vert z + \widehat{A}B^{-1}y\Vert_2\geq \Vert z + AB^{-1}y\Vert_2 - \Vert\widehat{A} - A\Vert_2\Vert B^{-1}\Vert_2\Vert y\Vert_2\geq\left(\frac{1}{C_{AB^{-1}}\sqrt{d}} - \frac{C_\delta\Vert B^{-1}\Vert_2}{\sqrt{n}}\right)\Vert y\Vert_2 - \Vert z\Vert_2 > M$, this implies
\begin{equation}
\sup_{\Vert z\Vert_\infty\leq x_0, \Vert y\Vert_2> n^{1/4}}\vert f_\epsilon(z + AB^{-1}y) - f_\epsilon(z + \widehat{A}B^{-1}y)\vert\leq 2\delta
\end{equation}
Combine with theorem \ref{thm1} and \eqref{Decop1} we prove the first result.

Notice that for any $\delta > 0$, there exists $M > 0$ such that $\int_{\Vert w\Vert_\infty\geq M \cup \Vert z\Vert_\infty \geq M}f_\eta(w)f_\eta(z)dwdz < \delta$, then
\begin{equation}
\begin{aligned}
\sup_{0<x\leq x_0}\vert G(x) - Prob\left(\Vert Y_{n+1} - BAB^{-1}Y_n\Vert_\infty\leq x\right)\vert\\
\leq \frac{1}{\vert\det(B)\vert}\sup_{0<x\leq x_0, w,z\in\mathbf{R}^d}\int_{\Vert y\Vert_\infty\leq x}\vert \widehat{f}_\epsilon(B^{-1}y + \widehat{A}B^{-1}w - B^{-1}z) - f_\epsilon(B^{-1}y + \widehat{A}B^{-1}w - B^{-1}z)\vert dy\\
+2\int_{\Vert w\Vert_\infty \geq M\cup\Vert z\Vert_\infty\geq M} f_\eta(w)f_\eta(z)dwdz\\
+ \frac{1}{\vert\det(B)\vert}\sup_{0<x\leq x_0, \Vert w\Vert_\infty,\Vert z\Vert_\infty\leq M}\int_{\Vert y\Vert_\infty\leq x}\vert f_\epsilon(B^{-1}y + \widehat{A}B^{-1}w - B^{-1}z) - f_\epsilon(B^{-1}y + AB^{-1}w - B^{-1}z)\vert dy\\
\leq \frac{\Vert \widehat{f}_\epsilon - f_\epsilon\Vert_{L_2}\times (2x_0)^{d/2}}{\sqrt{\vert\det(B)\vert}} + 2\delta \\ +\frac{(2x_0)^d}{\vert\det(B)\vert}\times \sup_{\Vert y\Vert_\infty\leq x_0, \Vert w\Vert_\infty,\Vert z\Vert_\infty\leq M}\vert f_\epsilon(B^{-1}y + \widehat{A}B^{-1}w - B^{-1}z) - f_\epsilon(B^{-1}y + AB^{-1}w - B^{-1}z)\vert
\end{aligned}
\end{equation}
If $\Vert\widehat{A} - A\Vert_2\leq C_\delta / \sqrt{n}$, then from remark \ref{remark1}, for sufficiently large $n$,
\begin{equation}
\begin{aligned}
\sup_{\Vert y\Vert_\infty\leq x_0, \Vert w\Vert_\infty,\Vert z\Vert_\infty\leq M}\vert f_\epsilon(B^{-1}y + \widehat{A}B^{-1}w - B^{-1}z) - f_\epsilon(B^{-1}y + AB^{-1}w - B^{-1}z)\vert\\
\leq \sup_{\Vert x - y\Vert_2\leq \frac{C_\delta\sqrt{d}M\Vert B^{-1}\Vert_2}{\sqrt{n}}}\vert f_\epsilon(x) - f_\epsilon(y)\vert < \delta
\end{aligned}
\end{equation}
which implies the second result.
\end{proof}

\begin{proof}[Proof of corollary \ref{MainCoro}]
Lemma \ref{lemma2} and remark \ref{remark2} imply the first result.

Notice that $Y_n = O_p(1)$, from theorem \ref{thm1} we know that $\Vert(\widehat{A} - A)B^{-1}Y_n\Vert_2\leq \Vert\widehat{A} - A\Vert_2\Vert B^{-1}\Vert_2\Vert Y_n\Vert_2 = O_p(1/\sqrt{n})$ and $\Vert B(\widehat{A} - A)B^{-1}Y_n\Vert_2\leq \Vert B\Vert_2\Vert\widehat{A} - A\Vert_2\Vert B^{-1}\Vert_2\Vert Y_n\Vert_2 = O_p(1/\sqrt{n})$, so for any $\delta > 0$, there exists a constant $C_\delta$ such that $Prob\left(\Vert(\widehat{A} - A)B^{-1}Y_n\Vert_2\leq \frac{C_\delta}{\sqrt{n}}\right) > 1-\delta$ and $Prob\left(\Vert B(\widehat{A} - A)B^{-1}Y_n\Vert_2\leq \frac{C_\delta}{\sqrt{n}}\right) > 1-\delta$ for any $n$. We adopt (1.2) in \cite{subsample} for and define
\begin{equation}
\begin{aligned}
c_{\mathcal{N}, 1-\alpha} = \inf\{x\geq 0| Prob\left(\Vert X_{2} - AB^{-1}Y_1\Vert_\infty\leq x\right)\geq 1-\alpha\}\\
c_{\mathcal{G}, 1-\alpha} = \inf\{x\geq 0| Prob\left(\Vert Y_{2} - BAB^{-1}Y_1\Vert_\infty\leq x\right)\geq 1-\alpha\}
\end{aligned}
\end{equation}
for all $0 <\alpha < 1$. From lemma \ref{lemma2} and dominated convergence theorem we know that $Prob\left(\Vert X_{2} - AB^{-1}Y_1\Vert_\infty\leq x\right)$ and $Prob\left(\Vert Y_{2} - BAB^{-1}Y_1\Vert_\infty\leq x\right)$ are continuous in $x > 0$, so $Prob\left(\Vert X_{2} - AB^{-1}Y_1\Vert_\infty\leq c_{\mathcal{N}, 1-\alpha}\right) = 1 - \alpha$ and $Prob\left(\Vert Y_{2} - BAB^{-1}Y_1\Vert_\infty\leq c_{\mathcal{G}, 1-\alpha}\right) = 1-\alpha$.

For a given $0<\alpha<1$, by choosing sufficiently large $x_0$ such that \eqref{Cond} happens and adopting notations in remark \ref{remark2}, for any $0 < \delta < \alpha / 8$ such that $1-\alpha - 2\delta > 0$, \eqref{ImportEq} implies
\begin{equation}
\begin{aligned}
N(x_0)\geq Prob\left(\Vert X_{n+1} - AB^{-1}Y_n\Vert_\infty\leq x_0\right) - \delta\geq 1-\frac{3}{8}\alpha\Rightarrow 0 < x^*_{\mathcal{N}}\leq x_0\\
Prob\left(\Vert X_{n+1} - AB^{-1}Y_n\Vert_\infty\leq x^*_{\mathcal{N}}\right)\geq 1-\alpha - \delta\Rightarrow x^*_{\mathcal{N}}\geq c_{\mathcal{N}, 1 - \alpha - \delta}\\
Prob\left(\Vert X_{n+1} - AB^{-1}Y_n\Vert_\infty\leq x^*_{\mathcal{N}}\right)\leq 1 - \alpha + \delta \Rightarrow x^*_{\mathcal{N}}\leq c_{\mathcal{N}, 1-\alpha +\delta}
\end{aligned}
\label{HHs}
\end{equation}
In \eqref{HHs} we consider $Prob\left(\Vert X_{n+1} - AB^{-1}Y_n\Vert_\infty\leq x\right)$ as a function in $x$ and plug in $x = x^*_{\mathcal{N}}$ instead of taking expectations. Similarly we can derive $x^*_{\mathcal{G}}\geq c_{\mathcal{G}, 1 - \alpha - \delta}$ and $x^*_{\mathcal{G}}\leq c_{\mathcal{G}, 1-\alpha +\delta}$. From theorem \ref{thm1}, we have for sufficiently large $n$
\begin{equation}
\begin{aligned}
Prob\left(\Vert X_{n+1} - \widehat{A}B^{-1}Y_n\Vert_\infty\leq x_{\mathcal{N}}^*\right)\\
\leq Prob\left(\Vert(\widehat{A} - A)B^{-1}Y_n\Vert_2 >  \frac{C_\delta}{\sqrt{n}}\right) + Prob\left(\sup_{0 < x\leq x_0}\vert N(x) - Prob\left(\Vert X_{n+1} - AB^{-1}Y_n\Vert_\infty\leq x\right)\vert > \delta\right)\\
+ Prob\left(\Vert X_{n+1} - AB^{-1}Y_n\Vert_\infty\leq c_{\mathcal{N}, 1-\alpha + \delta} + \frac{C_\delta}{\sqrt{n}}\right)\leq 1 - \alpha + 4\delta\\
Prob\left(\Vert X_{n+1} - \widehat{A}B^{-1}Y_n\Vert_\infty\leq x_{\mathcal{N}}^*\right)\\
\geq Prob\left(\Vert X_{n+1} - AB^{-1}Y_n\Vert_\infty\leq x_{\mathcal{N}}^* - \frac{C_\delta}{\sqrt{n}}\cap \Vert(\widehat{A} - A)B^{-1}Y_n\Vert_2 \leq  \frac{C_\delta}{\sqrt{n}}\right)\\
\geq Prob\left(\Vert X_{n+1} - AB^{-1}Y_n\Vert_\infty\leq c_{\mathcal{N}, 1-\alpha-\delta} - \frac{C_\delta}{\sqrt{n}}\cap\sup_{0 < x\leq x_0}\vert N(x) - Prob\left(\Vert X_{n+1} - AB^{-1}Y_n\Vert_\infty\leq x\right)\vert\leq \delta\right) - \delta\\
\geq 1-\alpha-4\delta
\end{aligned}
\label{Zels}
\end{equation}
which implies \eqref{ResCoro} by choosing $\delta \to 0$. Similarly we can replace $X_{n+1} - \widehat{A}B^{-1}Y_n$ and $x_{\mathcal{N}}^*$ in \eqref{Zels} by $Y_{n+1} - B\widehat{A}B^{-1}Y_n$ and $x^*_{\mathcal{G}}$ and derive the third result in \eqref{ResCoro}.
\end{proof}
\end{document}